# STABILITY RESULTS FOR SYSTEMS DESCRIBED BY RETARDED FUNCTIONAL DIFFERENTIAL EQUATIONS

**Iasson Karafyllis**[*], **Pierdomenico Pepe**[**] **and Zhong-Ping Jiang**[***]

[*] **Department of Environmental Engineering, Technical University of Crete, 73100, Chania, Greece**
email: ikarafyl@enveng.tuc.gr

[**]**Dipartimento di Ingegneria Elettrica, Universita degli Studi dell'Aquila, Monteluco di Roio, 67040, L'Aquila, Italy**
email: pepe@ing.univaq.it

[***]**Department of Electrical and Computer Engineering, Polytechnic University, Six Metrotech Center, Brooklyn, NY 11201, U.S.A.**
email: zjiang@control.poly.edu

**Abstract**
In this work characterizations of notions of output stability for uncertain time-varying systems described by retarded functional differential equations are provided. Particularly, characterizations by means of Lyapunov and Razumikhin functions of uniform and non-uniform in time Robust Global Asymptotic Output Stability and Input-to-Output Stability are given. The results of this work have been developed for systems with outputs in abstract normed linear spaces in order to allow outputs with no delay, with discrete or distributed delay or functional outputs with memory.

## 1. Introduction

The introduction of the notion of Input-to-State Stability (ISS) in [34] for finite-dimensional systems described by ordinary differential equations, led to an exceptionally rich period of progress in mathematical systems and control theory. The notion of ISS and its characterizations given in [35,36] were proved to be extremely useful for the expression of small-gain results (see [12,13,14,41,43]) and for the construction of robust feedback stabilizers (see for instance the textbooks [23,27,37]). The notion of ISS was extended to the notion of Input-to-Output Stability (IOS) in [39,40,10] and to the non-uniform in time notions of ISS and IOS in [18,19,21] (which extended the applicability of ISS to time-varying systems). Recently, semi-uniform notions of ISS have been proposed in [29]. The notions of ISS and IOS were recently proposed and characterized for discrete-time systems (see [15,16,20]) as well as to a wide class of systems with outputs (see [19]). It is our belief that the notions of ISS and IOS have become one of the most important conceptual tools for the development of nonlinear robust stability and control theory for a wide class of dynamical systems.

In this work we develop characterizations of various robust stability notions for uncertain systems described by Retarded Functional Differential Equations (RFDEs), including uniform and non-uniform in time ISS and IOS. The stability notions proposed in the present work are parallel to the robust stability notions used for finite-dimensional systems. Thus, it is expected that the results of the present paper will play an important role in mathematical systems and control theory for the important case of systems described by RFDEs.

Let $D \subseteq \Re^l$ be a non-empty set and $Y$ a normed linear space. We denote by $x(t)$ with $t \geq t_0$ the solution of the initial-value problem:

$$\dot{x}(t) = f(t, T_r(t)x, d(t)), \; t \geq t_0$$
$$Y(t) = H(t, T_r(t)x) \quad (1.1)$$
$$x(t) \in \Re^n, \; d(\cdot) \in M_D, \; Y(t) \in Y$$

with initial condition $T_r(t_0)x = x_0 \in C^0([-r,0]; \Re^n)$, where $r > 0$ is a constant $T_r(t)x$ denotes the "$r$-history" of $x$ at time $t$, i.e., $T_r(t)x := x(t+\theta); \theta \in [-r,0]$, and the mappings $f : \Re^+ \times C^0([-r,0]; \Re^n) \times D \to \Re^n$, $H : \Re^+ \times C^0([-r,0]; \Re^n) \to Y$ satisfy $f(t,0,d) = 0$, $H(t,0) = 0$ for all $(t,d) \in \Re^+ \times D$.

In this work we first study characterizations of notions of robust global asymptotic output stability for systems of the form (1.1) under weak assumptions (Section 3 of the paper). A major advantage of allowing the output to take values in abstract normed linear spaces is that using the framework of the case (1.1) we may consider:

- outputs with no delays, e.g. $Y(t) = h(t, x(t))$ with $Y = \Re^k$,

- outputs with discrete or distributed delay, e.g. $Y(t) = h(x(t), x(t-r))$ or $Y(t) = \int_{t-r}^{t} h(t, \theta, x(\theta)) d\theta$ with $Y = \Re^k$,

- functional outputs with memory, e.g. $Y(t) = h(t, \theta, x(t+\theta)); \theta \in [-r,0]$ or the identity output $Y(t) = T_r(t)x = x(t+\theta); \theta \in [-r,0]$ with $Y = C^0([-r,0]; \Re^k)$.

Moreover, it should be emphasized that our assumptions for (1.1) are very weak, since we do not assume boundedness or continuity of the right-hand side of the differential equation with respect to time or a Lipschitz condition for $f$. Furthermore, we do not assume that the disturbance set $D \subseteq \Re^l$ is compact.

Notions of output stability have been studied for finite-dimensional systems described by ordinary differential equations (see [39,40,10,19,21]). For systems described by RFDEs the notion of partial stability (which is a special case of the notion of global asymptotic output stability) has been studied in [2,3,9,44]. Particularly in [2], the authors provide Lyapunov characterizations of local partial stability for systems described by RFDEs without disturbances under the assumptions of the invariance of the attractive set and boundedness of the right-hand side of the differential equation with respect to time.

In this work we provide Lyapunov characterizations of Robust Global Asymptotic Output Stability (RGAOS) for systems described by RFDEs with disturbances (case (1.1)), without the hypothesis that the attractive set is invariant and without the assumption that the right-hand side of the differential equation (1.1) is bounded with respect to time. Paricularly, we consider uniform and non-uniform notions of RGAOS, which directly extend the corresponding notions of Robust Global Asymptotic Stability of an equilibrium point (see [3,5,8,22,24,25,26]).

We next continue with the study of the effect of non-vanishing perturbations for systems described by RFDEs (Section 4). Particularly, let $D \subseteq \Re^l$ be a non-empty set, $U \subseteq \Re^m$ a non-empty set with $0 \in U$ and $Y$ a normed linear space. We denote by $x(t)$ with $t \geq t_0$ the unique solution of the initial-value problem:

$$\dot{x}(t) = f(t, T_r(t)x, u(t), d(t)), \; t \geq t_0$$
$$Y(t) = H(t, T_r(t)x) \quad (1.2)$$
$$x(t) \in \Re^n, \; Y(t) \in Y, \; d(\cdot) \in M_D, \; u(\cdot) \in M_U$$

with initial condition $T_r(t_0)x = x_0 \in C^0([-r,0]; \Re^n)$, where $r > 0$ is a constant and the mappings $f : \Re^+ \times C^0([-r,0]; \Re^n) \times U \times D \to \Re^n$, $H : \Re^+ \times C^0([-r,0]; \Re^n) \to Y$ satisfy $f(t,0,0,d) = 0$, $H(t,0) = 0$ for all $(t,d) \in \Re^+ \times D$. In Section 4, we provide characterizations of the uniform and non-uniform in time IOS property for systems of the form (1.2) under weak hypotheses, which are usually satisfied in applications. The study of the uniform ISS property for autonomous systems described by RFDEs, was recently initiated in [33]. The technical issues that arise in the study of the case (1.2) are solved with a combination of the Lyapunov-like characterization given in [22] and an extension of the results in [32].

In Section 5, we develop Razumikhin conditions for the stability notions introduced in previous sections. The use of Razumikhin functions in the study of qualitative properties of the solutions of time-delay systems is emphasized in [8,31]. Recently, in [42] a major observation was established: Razumikhin theorems are "essentially" small-gain results. This idea is exploited in the present work, in order to produce novel results that are easily applicable.

Finally, in Section 6 we provide the concluding remarks of this work.

**Notations** Throughout this paper we adopt the following notations:

* Let $I \subseteq \Re$ be an interval. By $C^0(I;\Omega)$, we denote the class of continuous functions on $I$, which take values in $\Omega$. By $C^1(I;\Omega)$, we denote the class of functions on $I$ with continuous derivative, which take values in $\Omega$.

* For a vector $x \in \Re^n$ we denote by $|x|$ its usual Euclidean norm and by $x'$ its transpose. For $x \in C^0([-r,0];\Re^n)$ we define $\|x\|_r := \max_{\theta \in [-r,0]} |x(\theta)|$.

* $N$ denotes the set of positive integers and $\Re^+$ denotes the set of non-negative real numbers.

* We denote by $[R]$ the integer part of the real number $R$, i.e., the greatest integer, which is less than or equal to $R$.

* $\mathcal{E}$ denotes the class of non-negative $C^0$ functions $\mu: \Re^+ \to \Re^+$, for which it holds: $\int_0^{+\infty} \mu(t)dt < +\infty$ and $\lim_{t \to +\infty} \mu(t) = 0$.

* We denote by $K^+$ the class of positive $C^0$ functions defined on $\Re^+$. We say that a function $\rho: \Re^+ \to \Re^+$ is positive definite if $\rho(0) = 0$ and $\rho(s) > 0$ for all $s > 0$. By $K$ we denote the set of positive definite, increasing and continuous functions. We say that a positive definite, increasing and continuous function $\rho: \Re^+ \to \Re^+$ is of class $K_\infty$ if $\lim_{s \to +\infty} \rho(s) = +\infty$. By $KL$ we denote the set of all continuous functions $\sigma = \sigma(s,t): \Re^+ \times \Re^+ \to \Re^+$ with the properties: (i) for each $t \geq 0$ the mapping $\sigma(\cdot,t)$ is of class $K$; (ii) for each $s \geq 0$, the mapping $\sigma(s,\cdot)$ is non-increasing with $\lim_{t \to +\infty} \sigma(s,t) = 0$.

* Let $U \subseteq \Re^m$ be a non-empty set with $0 \in U$. By $B_U[0,r] := \{u \in U ; |u| \leq r\}$ we denote the closed sphere in $U \subseteq \Re^m$ with radius $r \geq 0$, centered at $0 \in U$.

* Let $D \subseteq \Re^l$ be a non-empty set. By $M_D$ we denote the class of all Lebesgue measurable and locally essentially bounded mappings $d: \Re^+ \to D$. By $\tilde{M}_D$ we denote the class of all right-continuous mappings $d: \Re^+ \to D$, with the property that there exists a countable set $A_d \subset \Re^+$ which is either finite or $A_d = \{t_k^d ; k = 1,...,\infty\}$ with $t_{k+1}^d > t_k^d > 0$ for all $k = 1,2,...$ and $\lim t_k^d = +\infty$, such that the mapping $t \in \Re^+ \setminus A_d \to d(t) \in D$ is continuous.

* Let $x:[a-r,b) \to \Re^n$ with $b > a > -\infty$ and $r > 0$. By $T_r(t)x$ we denote the "$r$-history" of $x$ at time $t \in [a,b)$, i.e., $T_r(t)x := x(t+\theta)$; $\theta \in [-r,0]$. Notice that $T_r(t)x \in C^0([-r,0];\Re^n)$.

* By $\| \ \|_Y$, we denote the norm of the normed linear space $Y$.

## 2. Main Assumptions and Preliminaries for Systems Described by RFDEs

In this work we consider systems of the form (1.1) under the following hypotheses:

**(H1)** The mapping $(x,d) \to f(t,x,d)$ is continuous for each fixed $t \geq 0$ and such that for every bounded $I \subseteq \Re^+$ and for every bounded $S \subset C^0([-r,0];\Re^n)$, there exists a constant $L \geq 0$ such that:

$$(x(0) - y(0))'(f(t,x,d) - f(t,y,d)) \leq L \max_{\tau \in [-r,0]} |x(\tau) - y(\tau)|^2 = L\|x-y\|_r^2$$
$$\forall t \in I, \forall (x,y) \in S \times S, \forall d \in D$$

Hypothesis (H1) is equivalent to the existence of a continuous function $L: \Re^+ \times \Re^+ \to \Re^+$ such that for each fixed $t \geq 0$ the mappings $L(t,\cdot)$ and $L(\cdot,t)$ are non-decreasing, with the following property:

$$(x(0)-y(0))'(f(t,x,d)-f(t,y,d)) \le L(t,\|x\|_r+\|y\|_r)\|x-y\|_r^2 \quad (2.1)$$
$$\forall (t,x,y,d) \in \Re^+ \times C^0([-r,0];\Re^n) \times C^0([-r,0];\Re^n) \times D$$

**(H2)** For every bounded $\Omega \subset \Re^+ \times C^0([-r,0];\Re^n)$ the image set $f(\Omega \times D) \subset \Re^n$ is bounded.

**(H3)** There exists a countable set $A \subset \Re^+$, which is either finite or $A = \{t_k\,; k = 1,...,\infty\}$ with $t_{k+1} > t_k > 0$ for all $k = 1,2,...$ and $\lim t_k = +\infty$, such that mapping $(t,x,d) \in (\Re^+ \setminus A) \times C^0([-r,0];\Re^n) \times D \to f(t,x,d)$ is continuous. Moreover, for each fixed $(t_0,x,d) \in \Re^+ \times C^0([-r,0];\Re^n) \times D$, we have $\lim_{t \to t_0^+} f(t,x,d) = f(t_0,x,d)$.

**(H4)** For every $\varepsilon > 0$, $t \in \Re^+$, there exists $\delta := \delta(\varepsilon,t) > 0$ such that $\sup\{|f(\tau,x,d)|\,; \tau \in \Re^+, d \in D, |\tau - t| + \|x\|_r < \delta\} < \varepsilon$.

**(H5)** The mapping $H(t,x)$ is locally Lipschitz, in the sense that for every bounded $I \subseteq \Re^+$ and for every bounded $S \subset C^0([-r,0];\Re^n)$, there exists a constant $L_H \ge 0$ such that:

$$\|H(t,x)-H(\tau,y)\|_Y \le L_H(|t-\tau|+\|x-y\|_r)$$
$$\forall (t,\tau) \in I \times I, \forall (x,y) \in S \times S$$

Hypothesis (H5) is equivalent to the existence of a continuous function $L_H : \Re^+ \times \Re^+ \to \Re^+$ such that for each fixed $t \ge 0$ the mappings $L_H(t,\cdot)$ and $L_H(\cdot,t)$ are non-decreasing, with the following property:

$$\|H(t,x)-H(\tau,y)\|_Y \le L_H\left(\max\{t,\tau\},\|x\|_r+\|y\|_r\right)(|t-\tau|+\|x-y\|_r) \quad (2.2)$$
$$\forall (t,\tau,x,y) \in \Re^+ \times \Re^+ \times C^0([-r,0];\Re^n) \times C^0([-r,0];\Re^n)$$

It is clear that (by virtue of hypotheses (H1-3) above and Lemma 1 in [7], page 4) for every $d \in M_D$ the composite map $f(t,x,d(t))$ satisfies the Caratheodory condition on $\Re^+ \times C^0([-r,0];\Re^n)$ and consequently, by virtue of Theorem 2.1 in [8] (and its extension given in paragraph 2.6 of the same book), for every $(t_0,x_0,d) \in \Re^+ \times C^0([-r,0];\Re^n) \times M_D$ there exists $h > 0$ and at least one continuous function $x : [t_0-r, t_0+h] \to \Re^n$, which is absolutely continuous on $[t_0, t_0+h]$ with $T_r(t_0)x = x_0$ and $\dot{x}(t) = f(t,T_r(t)x,d(t))$ almost everywhere on $[t_0, t_0+h]$. Let $x:[t_0-r,t_0+h] \to \Re^n$ and $y:[t_0-r,t_0+h] \to \Re^n$ two solutions of (1.1) with initial conditions $T_r(t_0)x = x_0$ and $T_r(t_0)y = y_0$ and corresponding to the same $d \in M_D$. Evaluating the derivative of the absolutely continuous map $z(t) = |x(t)-y(t)|^2$ on $[t_0, t_0+h]$ in conjunction with hypothesis (H1) above, we obtain the integral inequality:

$$|x(t)-y(t)|^2 \le |x(t_0)-y(t_0)|^2 + 2\int_{t_0}^t \tilde{L}\|T_r(\tau)x - T_r(\tau)y\|_r^2 d\tau, \quad \forall t \in [t_0, t_0+h]$$

where $\tilde{L} := L(t_0+h, a(t_0+h))$, $L(\cdot)$ is the function involved in (2.1) and $a(t) := \sup_{\tau \in [t_0-r,t]} |x(\tau)| + \sup_{\tau \in [t_0-r,t]} |y(\tau)|$. Consequently, we obtain:

$$\|T_r(t)(x-y)\|_r^2 \le \|x_0-y_0\|_r^2 + 2\int_{t_0}^t \tilde{L}\|T_r(\tau)(x-y)\|_r^2 d\tau, \quad \forall t \in [t_0, t_0+h]$$

and an immediate application of the Gronwall-Bellman inequality gives:

$$\|T_r(t)(x-y)\|_r \leq \|x_0 - y_0\|_r \exp\left(\tilde{L}(t-t_0)\right), \ \forall t \in [t_0, t_0+h] \quad (2.3)$$

Thus, we conclude that under hypotheses (H1-5), for every $(t_0, x_0, d) \in \Re^+ \times C^0([-r,0];\Re^n) \times M_D$ there exists $h > 0$ and exactly one continuous function $x:[t_0 - r, t_0 + h] \to \Re^n$, which is absolutely continuous on $[t_0, t_0+h]$ with $T_r(t_0)x = x_0$ and $\dot{x}(t) = f(t, T_r(t)x, d(t))$ almost everywhere on $[t_0, t_0+h]$. We denote by $\phi(t, t_0, x_0; d)$ the "$r$-history" of the unique solution of (1.1), i.e., $\phi(t, t_0, x_0; d) := T_r(t)x$, with initial condition $T_r(t_0)x = x_0$ corresponding to $d \in M_D$. Using hypothesis (H2) above and Theorem 3.2 in [8], we conclude that for every $(t_0, x_0, d) \in \Re^+ \times C^0([-r,0];\Re^n) \times M_D$ there exists $t_{\max} \in (t_0, +\infty]$, such that the unique solution $x(t)$ of (1.1) is defined on $[t_0 - r, t_{\max})$ and cannot be further continued. Moreover, if $t_{\max} < +\infty$ then we must necessarily have $\limsup_{t \to t_{\max}^-} |x(t)| = +\infty$. An immediate consequence of inequalities (2.3) and (2.2) is the following inequality which holds for every pair $\phi(\cdot, t_0, x_0; d):[t_0, t_{\max}^x) \to C^0([-r,0];\Re^n)$, $\phi(\cdot, t_0, y_0; d):[t_0, t_{\max}^y) \to C^0([-r,0];\Re^n)$ of solutions of (1.1) with initial conditions $T_r(t_0)x = x_0$, $T_r(t_0)y = y_0$, corresponding to the same $d \in M_D$ and for all $t \in [t_0, t_1)$ with $t_1 = \min\{t_{\max}^x; t_{\max}^y\}$:

$$\begin{aligned}
&\|\phi(t, t_0, x_0; d) - \phi(t, t_0, y_0; d)\|_r \leq \|x_0 - y_0\|_r \exp(L(t, a(t))(t-t_0)) \\
&\|H(t, \phi(t, t_0, x_0; d)) - H(t, \phi(t, t_0, x_0; d))\|_Y \leq L_H(t, a(t))\|x_0 - y_0\|_r \exp(L(t, a(t))(t-t_0)) \\
&a(t) = \sup_{\tau \in [t_0, t]} \left( \|\phi(t, t_0, x_0; d)\|_r + \|\phi(t, t_0, y_0; d)\|_r \right)
\end{aligned} \quad (2.4)$$

Since $f(t, 0, d) = 0$ for all $(t, d) \in \Re^+ \times D$, it follows that $\phi(t, t_0, 0; d) = 0 \in C^0([-r,0];\Re^n)$ for all $(t_0, d) \in \Re^+ \times M_D$ and $t \geq t_0$. Furthermore, (2.4) implies that for every $\varepsilon > 0$, $T, h \geq 0$ there exists $\delta := \delta(\varepsilon, T, h) > 0$ such that:

$$\|x\|_r < \delta \quad \Rightarrow \quad \sup\{\|\phi(\tau, t_0, x; d)\|_r \ ; d \in M_D, \tau \in [t_0, t_0+h], t_0 \in [0,T]\} < \varepsilon$$

Thus $0 \in C^0([-r,0];\Re^n)$ is a robust equilibrium point for (1.1) in the sense described in [19].

It should be emphasized that if $d \in \tilde{M}_D$ then the map $t \to f(t, x, d(t))$ is right-continuous on $\Re^+$ and continuous on $\Re^+ \setminus (A \cup A_d)$. Applying repeatedly Theorem 2.1 in [8] on each one of the intervals contained in $[t_0, t_{\max}) \setminus (A \cup A_d)$, we conclude that the solution satisfies $\dot{x}(t) = f(t, T_r(t)x, d(t))$ for all $t \in [t_0, t_{\max}) \setminus (A \cup A_d)$. By virtue of the mean value theorem, it follows that $\lim_{h \to 0^+} \frac{x(t+h) - x(t)}{h} = f(t, T_r(t)x, d(t))$ for all $t \in [t_0, t_{\max})$.

An important property for systems of the form (1.1) is Robust Forward Completeness (RFC) (see [19]). This property will be used extensively in the following sections of the present work.

**Definition 2.1:** *We say that (1.1) under hypotheses (H1-5) is **Robustly Forward Complete (RFC)** if for every $s \geq 0$, $T \geq 0$, it holds that*

$$\sup\{\|\phi(t_0+\xi, t_0, x_0; d)\|_r \ ; \ \xi \in [0,T], \|x_0\|_r \leq s, t_0 \in [0,T], d \in M_D \} < +\infty$$

Concerning systems of the form (1.2) the following hypotheses will be valid throughout the text:

**(S1)** The mapping $(x, u, d) \to f(t, x, u, d)$ is continuous for each fixed $t \geq 0$ and such that for every bounded $I \subseteq \Re^+$ and for every bounded $S \subset C^0([-r,0];\Re^n) \times U$, there exists a constant $L \geq 0$ such that:

$$(x(0) - y(0))'(f(t,x,u,d) - f(t,y,u,d)) \leq L \max_{\tau \in [-r,0]} |x(\tau) - y(\tau)|^2 = L \|x - y\|_r^2$$
$$\forall t \in I, \forall (x,u,y,u) \in S \times S, \forall d \in D$$

Hypothesis (S1) is equivalent to the existence of a continuous function $L : \Re^+ \times \Re^+ \to \Re^+$ such that for each fixed $t \geq 0$ the mappings $L(t,\cdot)$ and $L(\cdot,t)$ are non-decreasing, with the following property:

$$(x(0) - y(0))'(f(t,x,u,d) - f(t,y,u,d)) \leq L(t, \|x\|_r + \|y\|_r + |u|) \|x - y\|_r^2$$
$$\forall (t,x,y,d,u) \in \Re^+ \times C^0([-r,0];\Re^n) \times C^0([-r,0];\Re^n) \times D \times U \quad (2.5)$$

**(S2)** For every bounded $\Omega \subset \Re^+ \times C^0([-r,0];\Re^n) \times U$ the image set $f(\Omega \times D) \subset \Re^n$ is bounded.

**(S3)** There exists a countable set $A \subset \Re^+$, which is either finite or $A = \{t_k \,; k = 1,...,\infty\}$ with $t_{k+1} > t_k > 0$ for all $k = 1,2,...$ and $\lim t_k = +\infty$, such that mapping $(t,x,u,d) \in (\Re^+ \setminus A) \times C^0([-r,0];\Re^n) \times U \times D \to f(t,x,u,d)$ is continuous. Moreover, for each fixed $(t_0,x,u,d) \in \Re^+ \times C^0([-r,0];\Re^n) \times U \times D$, we have $\lim_{t \to t_0^+} f(t,x,u,d) = f(t_0,x,u,d)$.

**(S4)** For every $\varepsilon > 0$, $t \in \Re^+$, there exists $\delta := \delta(\varepsilon, t) > 0$ such that $\sup\{|f(\tau,x,u,d)| \,; \tau \in \Re^+, d \in D, u \in U, |\tau - t| + \|x\|_r + |u| < \delta\} < \varepsilon$.

**(S5)** The mapping $u \to f(t,x,u,d)$ is locally Lipschitz, in the sense that for every bounded $I \subseteq \Re^+$ and for every bounded $S \subset C^0([-r,0];\Re^n) \times U$, there exists a constant $L_U \geq 0$ such that:

$$|f(t,x,u,d) - f(t,x,v,d)| \leq L_U |u-v|$$
$$\forall t \in I, \forall (x,u,x,v) \in S \times S, \forall d \in D$$

Hypothesis (S5) is equivalent to the existence of a continuous function $L_U : \Re^+ \times \Re^+ \to \Re^+$ such that for each fixed $t \geq 0$ the mappings $L_U(t,\cdot)$ and $L_U(\cdot,t)$ are non-decreasing, with the following property:

$$|f(t,x,u,d) - f(t,x,v,d)| \leq L_U(t, \|x\|_r + |u| + |v|) |u-v|$$
$$\forall (t,x,d,u,v) \in \Re^+ \times C^0([-r,0];\Re^n) \times D \times U \times U \quad (2.6)$$

**(S6):** $U$ is a positive cone, i.e., for all $u \in U$ and $\lambda \geq 0$ it follows that $(\lambda u) \in U$.

**(S7)** The mapping $H(t,x)$ satisfies hypothesis (H5) given above.

It is clear that (by virtue of hypotheses (S1-3) above and Lemma 1 in [7], page 4) for every $(d,u) \in M_D \times M_U$ the composite map $f(t,x,u(t),d(t))$ satisfies the Caratheodory condition on $\Re^+ \times C^0([-r,0];\Re^n)$ and consequently, by virtue of Theorem 2.1 in [8] (and its extension given in paragraph 2.6 of the same book), for every $(t_0,x_0,u,d) \in \Re^+ \times C^0([-r,0];\Re^n) \times M_U \times M_D$ there exists $h > 0$ and at least one continuous function $x : [t_0 - r, t_0 + h] \to \Re^n$, which is absolutely continuous on $[t_0, t_0 + h]$ with $T_r(t_0)x = x_0$ and $\dot{x}(t) = f(t, T_r(t)x, u(t), d(t))$ almost everywhere on $[t_0, t_0 + h]$. Let $x : [t_0 - r, t_0 + h] \to \Re^n$ and $y : [t_0 - r, t_0 + h] \to \Re^n$ two solutions of (1.2) with initial conditions $T_r(t_0)x = x_0$ and $T_r(t_0)y = y_0$ and corresponding to the same $(d,u) \in M_D \times M_U$. Evaluating the derivative of the absolutely continuous map $z(t) = |x(t) - y(t)|^2$ on $[t_0, t_0 + h]$ in conjunction with hypothesis (H1) above, we obtain the integral inequality:

$$|x(t) - y(t)|^2 \le |x(t_0) - y(t_0)|^2 + 2\int_{t_0}^{t} \tilde{L} \|T_r(\tau)x - T_r(\tau)y\|_r^2 d\tau, \ \forall t \in [t_0, t_0 + h]$$

where $\tilde{L} := L(t_0 + h, a(t_0 + h))$, $L(\cdot)$ is the function involved in (2.5) and $a(t) := \sup_{\tau \in [t_0-r,t]} |x(\tau)| + \sup_{\tau \in [t_0-r,t]} |y(\tau)| + \sup_{\tau \in [t_0,t]} |u(\tau)|$. Consequently, we obtain:

$$\|T_r(t)(x-y)\|_r^2 \le \|x_0 - y_0\|_r^2 + 2\int_{t_0}^{t} \tilde{L} \|T_r(\tau)(x-y)\|_r^2 d\tau, \ \forall t \in [t_0, t_0 + h]$$

and immediate application of the Gronwall-Bellman inequality gives:

$$\|T_r(t)(x-y)\|_r \le \|x_0 - y_0\|_r \exp(\tilde{L}(t-t_0)), \ \forall t \in [t_0, t_0 + h] \quad (2.7)$$

Thus, we conclude that under hypotheses (S1-7), for every $(t_0, x_0, d, u) \in \Re^+ \times C^0([-r,0]; \Re^n) \times M_D \times M_U$ there exists $h > 0$ and exactly one continuous function $x : [t_0 - r, t_0 + h] \to \Re^n$, which is absolutely continuous on $[t_0, t_0 + h]$ with $T_r(t_0)x = x_0$ and $\dot{x}(t) = f(t, T_r(t)x, u(t), d(t))$ almost everywhere on $[t_0, t_0 + h]$. We denote by $\phi(t, t_0, x_0; u, d)$ the "$r$-history" of the unique solution of (1.2), i.e., $\phi(t, t_0, x_0; u, d) := T_r(t)x$, with initial condition $T_r(t_0)x = x_0$ corresponding to $(d, u) \in M_D \times M_U$. Using hypothesis (S2) above and Theorem 3.2 in [8], we conclude that for every $(t_0, x_0, d, u) \in \Re^+ \times C^0([-r,0]; \Re^n) \times M_D \times M_U$ there exists $t_{\max} \in (t_0, +\infty]$, such that the unique solution $x(t)$ of (1.2) is defined on $[t_0 - r, t_{\max})$ and cannot be further continued. Moreover, if $t_{\max} < +\infty$ then we must necessarily have $\limsup_{t \to t_{\max}} |x(t)| = +\infty$. An immediate consequence of inequalities (2.7) and (2.2) is the following inequality, which holds for every pair $\phi(\cdot, t_0, x_0; u, d) : [t_0, t_{\max}^x) \to C^0([-r,0]; \Re^n)$, $\phi(\cdot, t_0, y_0; u, d) : [t_0, t_{\max}^y) \to C^0([-r,0]; \Re^n)$ of solutions of (1.2) with initial conditions $T_r(t_0)x = x_0$, $T_r(t_0)y = y_0$, corresponding to the same $(d, u) \in M_D \times M_U$ and for all $t \in [t_0, t_1)$ with $t_1 = \min\{t_{\max}^x; t_{\max}^y\}$:

$$\|\phi(t, t_0, x_0; u, d) - \phi(t, t_0, y_0; u, d)\|_r \le \|x_0 - y_0\|_r \exp(L(t, a(t))(t - t_0))$$
$$\|H(t, \phi(t, t_0, x_0; u, d)) - H(t, \phi(t, t_0, x_0; u, d))\|_Y \le L_H(t, a(t))\|x_0 - y_0\|_r \exp(L(t, a(t))(t - t_0)) \quad (2.8)$$
$$a(t) = \sup_{\tau \in [t_0,t]} \left( \|\phi(t, t_0, x_0; u, d)\|_r + \|\phi(t, t_0, y_0; u, d)\|_r \right) + \sup_{\tau \in [t_0,t]} |u(\tau)|$$

It should be emphasized that if $(d, u) \in \tilde{M}_D \times \tilde{M}_U$ then the map $t \to f(t, x, u(t), d(t))$ is right-continuous on $\Re^+$ and continuous on $\Re^+ \setminus (A \cup A_d \cup A_u)$. Applying repeatedly Theorem 2.1 in [8] on each one of the intervals contained in $[t_0, t_{\max}) \setminus (A \cup A_d \cup A_u)$, we conclude that the solution satisfies $\dot{x}(t) = f(t, T_r(t)x, u(t), d(t))$ for all $t \in [t_0, t_{\max}) \setminus (A \cup A_d \cup A_u)$. By virtue of the mean value theorem, it follows that $\lim_{h \to 0^+} \frac{x(t+h) - x(t)}{h} = f(t, T_r(t)x, u(t), d(t))$ for all $t \in [t_0, t_{\max})$.

An important remark concerning hypotheses (S1-7) is that for the case $u \equiv 0$ we obtain a system of the form (1.1) which satisfies hypotheses (H1-5). The same conclusion holds if $u = k(t, T_r(t)x)$, where $k : \Re^+ \times C^0([-r,0]; \Re^n) \to U$ is a mapping which is Lipschitz on bounded sets with $k(t,0) = 0$ for all $t \ge 0$.

An important property for systems of the form (1.2) is Robust Forward Completeness (RFC) from an external input (see [19]). This property will be used extensively in the following sections of the present work. Notice that the notion of Robust Forward Completeness (RFC) from the input $u \in M_U$ coincides with the notion of Robust Forward Completeness (RFC) for systems of the form (1.2) when $u \equiv 0$.

**Definition 2.2:** *We say that (1.2) under hypotheses (S1-7) is **robustly forward complete (RFC) from the input** $u \in M_U$ if for every $s \geq 0$, $T \geq 0$, it holds that*

$$\sup\left\{\|\phi(t_0 + \xi, t_0, x_0; u, d)\|_r ; u \in M_{B_U[0,s]}, \xi \in [0,T], \|x_0\|_r \leq s, t_0 \in [0,T], d \in M_D \right\} < +\infty$$

In order to study the asymptotic properties of the solutions of systems of the form (1.1) or (1.2), we will use Lyapunov functionals and functions. Therefore, a detailed list of certain notions and properties concerning functionals is needed.

Let $x \in C^0\left([-r,0]; \Re^n\right)$. By $E_h(x; v)$, where $0 \leq h < r$ and $v \in \Re^n$ we denote the following operator:

$$E_h(x; v) := \begin{cases} x(0) + (\theta + h) v & \text{for } -h < \theta \leq 0 \\ x(\theta + h) & \text{for } -r \leq \theta \leq -h \end{cases} \quad (2.9)$$

Let $V : \Re^+ \times C^0\left([-r,0]; \Re^n\right) \to \Re$. We define

$$V^0(t, x; v) := \limsup_{\substack{h \to 0^+ \\ y \to 0, y \in C^0([-r,0]; \Re^n)}} \frac{V(t + h, E_h(x; v) + hy) - V(t, x)}{h} \quad (2.10)$$

**Remark 2.3:** For mappings $V : \Re^+ \times C^0\left([-r,0]; \Re^n\right) \to \Re$, which are Lipschitz on bounded sets of $\Re^+ \times C^0\left([-r,0]; \Re^n\right)$, the above derivative coincides with the derivative introduced in [6] and was used later in [4].

The following lemma presents some elementary properties of the generalized derivative given above. Notice that the function $(t, x, v) \to V^0(t, x; v)$ may take values in the extended real number set $\Re^* = [-\infty, +\infty]$. Its proof is almost identical with Lemma 2.7 in [22]. Notice that we are not assuming that the mapping $V : \Re^+ \times C^0\left([-r,0]; \Re^n\right) \to \Re$ is Lipschitz on bounded sets of $\Re^+ \times C^0\left([-r,0]; \Re^n\right)$.

**Lemma 2.4:** *Let* $V : \Re^+ \times C^0\left([-r,0]; \Re^n\right) \to \Re$ *and let* $x \in C^0([t_0 - r, t_{\max}); \Re^n)$ *a solution of (1.1) under hypotheses (H1-5) corresponding to certain* $d \in M_D$ *(or a solution of (1.2) under hypotheses (S1-7) corresponding to certain* $(d, u) \in M_D \times M_U$), *where* $t_{\max} \in (t_0, +\infty]$ *is the maximal existence time of the solution. Then it holds that*

$$\limsup_{h \to 0^+} h^{-1}\left(V(t + h, T_r(t + h)x) - V(t, T_r(t)x)\right) \leq V^0(t, T_r(t)x; D^+ x(t)), \text{ a.e. on } [t_0, t_{\max}) \quad (2.11)$$

*where* $D^+ x(t) = \lim_{h \to 0^+} h^{-1}\left(x(t + h) - x(t)\right)$. *Moreover, if* $d \in \tilde{M}_D$ *(or* $(d, u) \in \tilde{M}_D \times \tilde{M}_U$) *then (2.11) holds for all* $t \in [t_0, t_{\max})$.

**Proof** It suffices to show that (2.11) holds for all $t \in [t_0, t_{\max}) \setminus I$ where $I \subset [t_0, t_{\max})$ is the set of zero Lebesgue measure such that $D^+ x(t) = \lim_{h \to 0^+} h^{-1}\left(x(t + h) - x(t)\right)$ is not defined on $I$. Let $h > 0$ and $t \in [t_0, t_{\max}) \setminus I$. We define:

$$T_r(t + h)x - E_h(T_r(t)x; D^+ x(t)) = h y_h \quad (2.12)$$

where

$$y_h = h^{-1}\begin{cases} x(t + h + \theta) - x(t) - (\theta + h)D^+ x(t) & \text{for } -h < \theta \leq 0 \\ 0 & \text{for } -r \leq \theta \leq -h \end{cases}$$

and notice that $y_h \in C^0([-r,0]; \Re^n)$ (as difference of continuous functions, see (2.12) above). Equivalently $y_h$ satisfies:

$$y_h := \begin{cases} \dfrac{\theta+h}{h}\left(\dfrac{x(t+\theta+h)-x(t)}{\theta+h} - D^+x(t)\right) & \text{for } -h < \theta \leq 0 \\ 0 & \text{for } -r \leq \theta \leq -h \end{cases}$$

with $\|y_h\|_r \leq \sup\left\{\left|\dfrac{x(t+s)-x(t)}{s} - D^+x(t)\right|; 0 < s \leq h\right\}$. Since $\lim_{h\to 0^+}\dfrac{x(t+h)-x(t)}{h} = D^+x(t)$ we obtain that $y_h \to 0$ as $h \to 0^+$. Finally, by virtue of definitions (2.10), (2.12) and since $y_h \to 0$ as $h \to 0^+$, we have:

$$\limsup_{h\to 0^+} h^{-1}\left(V(t+h, T_r(t+h)x) - V(t, T_r(t)x)\right)$$
$$= \limsup_{h\to 0^+} h^{-1}\left(V(t+h, E_h(T_r(t)x; D^+x(t)) + hy_h) - V(t, T_r(t)x)\right) \leq V^0(t, T_r(t)x; D^+x(t))$$

The proof is complete. ◁

An important class of functionals is presented next.

**Definition 2.5:** *We say that a continuous functional $V : \Re^+ \times C^0([-r,0]; \Re^n) \to \Re^+$, is "almost Lipschitz on bounded sets", if there exist non-decreasing functions $M : \Re^+ \to \Re^+$, $P : \Re^+ \to \Re^+$, $G : \Re^+ \to [1,+\infty)$ such that for all $R \geq 0$, the following properties hold:*

**(P1)** *For every $x, y \in \{x \in C^0([-r,0]; \Re^n); \|x\|_r \leq R\}$, it holds that:*

$$|V(t,y) - V(t,x)| \leq M(R)\|y-x\|_r, \quad \forall t \in [0, R]$$

**(P2)** *For every absolutely continuous function $x : [-r,0] \to \Re^n$ with $\|x\|_r \leq R$ and essentially bounded derivative, it holds that:*

$$|V(t+h, x) - V(t,x)| \leq hP(R)\left(1 + \sup_{-r \leq \tau \leq 0}|\dot{x}(\tau)|\right), \text{ for all } t \in [0, R] \text{ and } 0 \leq h \leq \dfrac{1}{G\left(R + \sup_{-r \leq \tau \leq 0}|\dot{x}(\tau)|\right)}$$

For the important class of functionals which are almost Lipschitz on bounded sets we are in a position to prove a novel result, which extends the result of Theorem 4 in [32].

**Lemma 2.6:** *Let $V : \Re^+ \times C^0([-r,0]; \Re^n) \to \Re$ be a functional which is almost Lipschitz on bounded sets and let $x \in C^0([t_0-r, t_{\max}); \Re^n)$ a solution of (1.1) under hypotheses (H1-5) corresponding to certain $d \in M_D$ (or a solution of (1.2) under hypotheses (S1-7) corresponding to certain $(d,u) \in M_D \times M_U$) with initial condition $T_r(t_0)x = x_0 \in C^1([-r,0]; \Re^n)$, where $t_{\max} \in (t_0, +\infty]$ is the maximal existence time of the solution. Then for every $T \in (t_0, t_{\max})$, the mapping $[t_0, T] \ni t \to V(t, T_r(t)x)$ is absolutely continuous.*

**Proof:** It suffices to show that for every $T \in (t_0, t_{\max})$ and $\varepsilon > 0$ there exists $\delta > 0$ such that $\sum_{k=1}^{N}|V(b_k, T_r(b_k)x) - V(a_k, T_r(a_k)x)| < \varepsilon$ for every finite collection of pairwise disjoint intervals $[a_k, b_k] \subset [t_0, T]$ ($k = 1, ..., N$) with $\sum_{k=1}^{N}(b_k - a_k) < \delta$.

Let $T \in (t_0, t_{\max})$ and $\varepsilon > 0$ (arbitrary). Since the solution $x \in C^0([t_0-r, T]; \Re^n)$ of (1.1) under hypotheses (H1-5) corresponding to certain $d \in M_D$ (or the solution $x \in C^0([t_0-r, T]; \Re^n)$ of (1.2) under hypotheses (S1-7) corresponding to certain $(d,u) \in M_D \times M_U$) with initial condition $T_r(t_0)x = x_0 \in C^1([-r,0];\Re^n)$ is bounded on $[t_0 - r, T]$, there exists $R_1 > 0$ such that $\sup_{t_0 \leq \tau \leq T} \|T_r(\tau)x\|_r \leq R_1$. Moreover, by virtue of hypothesis (H2) (or hypothesis (S2)) and since $T_r(t_0)x = x_0 \in C^1([-r,0];\Re^n)$, there exists $R_2 > 0$ such that $\sup_{t_0 - r \leq \tau \leq T} |\dot{x}(\tau)| \leq R_2$. The previous observations in conjunction with properties (P1), (P2) of Definition 2.5 imply for every interval $[a,b] \subset [t_0, T]$ with $b - a \leq \dfrac{1}{G(R + R_2)}$:

$$|V(b, T_r(b)x) - V(a, T_r(a)x)| \leq (b-a) P(R_1)(1 + R_2) + M(R_1) \|T_r(b)x - T_r(a)x\|_r$$

In addition, the estimate $\sup_{t_0 - r \leq \tau \leq T} |\dot{x}(\tau)| \leq R_2$ implies $\|T_r(b)x - T_r(a)x\|_r \leq (b-a) R_2$ for every interval $[a,b] \subset [t_0, T]$. Consequently, we obtain for every interval $[a,b] \subset [t_0, T]$ with $b - a \leq \dfrac{1}{G(R + R_2)}$:

$$|V(b, T_r(b)x) - V(a, T_r(a)x)| \leq (b-a) \left[ P(R_1)(1 + R_2) + M(R_1) R_2 \right]$$

The previous inequality implies that for every finite collection of pairwise disjoint intervals $[a_k, b_k] \subset [t_0, T]$ ($k = 1,...,N$) with $\sum_{k=1}^{N}(b_k - a_k) < \delta$, where $\delta = \dfrac{1}{2} \min \left\{ \dfrac{1}{G(R + R_2)} ; \dfrac{\varepsilon}{P(R_1)(1 + R_2) + M(R_1) R_2} \right\} > 0$, it holds that $\sum_{k=1}^{N} |V(b_k, T_r(b_k)x) - V(a_k, T_r(a_k)x)| < \varepsilon$. The proof is complete. ◁

The following lemma extends the result presented in [32] and shows that appropriate estimates of the solutions of systems (1.1) and (1.2) hold globally. The proof of the following lemma is analogous to the proof of Proposition 2 in [32].

**Lemma 2.7:** *Suppose that there exist mappings* $\beta_1 : \Re^+ \times C^0([-r,0];\Re^n) \to \Re$, $\beta_2 : \Re^+ \times \Re^+ \times C^0([-r,0];\Re^n) \times A \to \Re$, *where* $A \subseteq M_D \times M_U$, *with the following properties:*

**(i)** *for each* $(t, t_0, d, u) \in \Re^+ \times \Re^+ \times A$, *the mappings* $x \to \beta_1(t,x)$, $x \to \beta_2(t, t_0, x, d, u)$ *are continuous,*

**(ii)** *there exists a continuous function* $M : \Re^+ \times \Re^+ \to \Re^+$ *such that*

$$\sup \left\{ \beta_2(t_0 + \xi, t_0, x_0, d, u) ; \sup_{t \geq 0} |u(\tau)| \leq s, \xi \in [0,T], x_0 \in C^0([-r,0];\Re^n), \|x_0\|_r \leq s, t_0 \in [0,T], (d,u) \in A \right\} \leq M(T,s)$$

**(iii)** *for every* $(t_0, x_0, d, u) \in \Re^+ \times C^1([-r,0];\Re^n) \times A$ *the solution* $x(t)$ *of (1.2) with initial condition* $T_r(t_0)x = x_0$ *corresponding to input* $(d,u) \in A$ *satisfies:*

$$\beta_1\big(t, T_r(t)x\big) \leq \beta_2\big(t, t_0, x_0, d, u\big), \ \forall t \geq t_0 \tag{2.13}$$

*Moreover, suppose that one of the following properties holds:*

**(iv)** $c(T,s) := \sup \left\{ \|T_r(t_0 + \xi)x\|_r ; \sup_{t \geq 0} |u(\tau)| \leq s, \xi \in [0,T], x_0 \in C^0([-r,0];\Re^n), \|x_0\|_r \leq s, t_0 \in [0,T], (d,u) \in A \right\} < +\infty$

**(v)** there exist functions $a \in K_\infty$, $\mu \in K^+$ and a constant $R \geq 0$ such that $a(\mu(t)|x(0)|) \leq \beta_1(t, x) + R$ for all $(t, x) \in \Re^+ \times C^0([-r,0]; \Re^n)$

*Then for every* $(t_0, x_0, d, u) \in \Re^+ \times C^0([-r,0]; \Re^n) \times A$ *the solution* $x(t)$ *of (1.2) with initial condition* $T_r(t_0)x = x_0$ *corresponding to input* $(d, u) \in A$ *exists for all* $t \geq t_0$ *and satisfies (2.13).*

**Remark 2.8:** Notice that the statement of Lemma 2.7 covers the case (1.1) since (1.1) can be considered as a system of the form (1.2) described by the equation $\dot{x}(t) = f(t, T_r(t)x, d(t)) + u(t)$, where $U = \Re^n$. Moreover, every solution of system (1.1) can be considered as a solution of system $\dot{x}(t) = f(t, T_r(t)x, d(t)) + u(t)$ with $(d, u) \in A := M_D \times \{0\}$.

**Proof of Lemma 2.7:** We distinguish the following cases:

(a) Property (iv) holds. The proof will be made by contradiction. Suppose on the contrary that there exists $(t_0, x_0, d, u) \in \Re^+ \times C^0([-r,0]; \Re^n) \times A$ and $t_1 > t_0$ such that the solution $x(t)$ of (1.2) with initial condition $T_r(t_0)x = x_0$ corresponding to input $(d, u) \in A$ satisfies:

$$\beta_1(t_1, T_r(t_1)x) > \beta_2(t_1, t_0, x_0, d, u)$$

Using (2.8) and property (iv) we obtain for all $\tilde{x}_0 \in C^0([-r,0]; \Re^n)$ with $\|x_0 - \tilde{x}_0\|_r \leq 1$:

$$\|T_r(t_1)x - T_r(t_1)\tilde{x}\|_r \leq \|x_0 - \tilde{x}_0\|_r \exp(L(t_1, \bar{c})(t_1 - t_0)) \tag{2.14}$$

where $\tilde{x}(t)$ denotes the solution of (1.2) with initial condition $T_r(t_0)x = \tilde{x}_0$ corresponding to input $(d, u) \in A$ and $\bar{c} = 2c\left(t_1, \|x_0\|_r + 1 + \sup_{t_0 \leq \tau \leq t_1} |u(\tau)|\right) + \sup_{t_0 \leq \tau \leq t_1} |u(\tau)|$. Notice that in order to obtain inequality (2.14) we have also used the causality argument that the solutions $x(t)$, $\tilde{x}(t)$ of (1.2) with initial condition $T_r(t_0)x = x_0$ and $T_r(t_0)x = \tilde{x}_0$, respectively, corresponding to input $(d, u) \in A$, depend only on the values of the input $u$ on the interval $[t_0, t]$.

Let $\varepsilon := \beta_1(t_1, T_r(t_1)x) - \beta_2(t_1, t_0, x_0, u) > 0$. Using property (iv), (2.14), density of $C^1([-r,0]; \Re^n)$ in $C^0([-r,0]; \Re^n)$, continuity of the mappings $x \to \beta_1(t_1, x)$, $x \to \beta_2(t_1, t_0, x, d, u)$, we conclude that there exists $\tilde{x}_0 \in C^1([-r,0]; \Re^n)$ such that:

$$\|x_0 - \tilde{x}_0\|_r \leq 1; \; |\beta_2(t_1, t_0, x_0, d, u) - \beta_2(t_1, t_0, \tilde{x}_0, d, u)| \leq \frac{\varepsilon}{2} \; ; \; |\beta_1(t_1, T_r(t_1)x) - \beta_1(t_1, T_r(t_1)\tilde{x})| \leq \frac{\varepsilon}{2}$$

where $\tilde{x}(t)$ denotes the solution of (1.2) with initial condition $T_r(t_0)x = \tilde{x}_0$ corresponding to input $(d, u) \in A$. Combining property (iii) for $\tilde{x}(t)$ with the above inequalities and the definition of $\varepsilon$ we obtain $\beta_1(t_1, T_r(t_1)x) > \beta_1(t_1, T_r(t_1)x)$, a contradiction.

(b) Property (v) holds. It suffices to show that property (iv) holds. Since there exist functions $a \in K_\infty$, $\mu \in K^+$ and a constant $R \geq 0$ such that $a(\mu(t)|x(0)|) \leq \beta_1(t, x) + R$ for all $(t, x) \in \Re^+ \times C^0([-r,0]; \Re^n)$, it follows that from property (iii) that for every $(t_0, x_0, d, u) \in \Re^+ \times C^1([-r,0]; \Re^n) \times A$ the solution $x(t)$ of (1.2) with initial condition $T_r(t_0)x = x_0$ corresponding to input $(d, u) \in A$ satisfies:

$$a(\mu(t)|x(t)|) \leq R + \beta_2(t, t_0, x_0, d, u), \; \forall t \geq t_0$$

Moreover, making use of property (ii), the above inequality, we obtain that for every $(t_0, x_0, d, u) \in \Re^+ \times C^1([-r,0]; \Re^n) \times A$ the solution $x(t)$ of (1.2) with initial condition $T_r(t_0)x = x_0$ corresponding to input $(d, u) \in A$ satisfies:

$$\left\| T_r(t)x \right\|_r \leq \left\| x_0 \right\|_r + 1 + \frac{1}{\mu(t)} a^{-1}\left( R + M\left( t, \left\| x_0 \right\|_r + \sup_{t_0 \leq \tau \leq t} \left| u(\tau) \right| \right) \right), \quad \forall t \geq t_0 \qquad (2.15)$$

Notice that in order to obtain inequality (2.15) we have also used the causality argument that the solution $x(t)$ of (1.2) with initial condition $T_r(t_0)x = x_0$ corresponding to input $(d,u) \in A$, depends only on the values of the input $u$ on the interval $[t_0, t]$.

We claim that estimate (2.15) holds for all $(t_0, x_0, d, u) \in \Re^+ \times C^0([-r,0]; \Re^n) \times A$. Notice that this claim implies directly that property (iv) holds with $c(T,s) := s + 1 + \dfrac{1}{\min_{0 \leq \tau \leq 2T} \mu(\tau)} a^{-1}\left( R + \max_{\substack{0 \leq x \leq 2s, \\ 0 \leq \tau \leq 2T}} M(\tau, s) \right)$. The proof of the claim will be made by contradiction. Suppose on the contrary that there exists $(t_0, x_0, d, u) \in \Re^+ \times C^0([-r,0]; \Re^n) \times A$ and $t_1 > t_0$ such that the solution $x(t)$ of (1.2) with initial condition $T_r(t_0)x = x_0$ corresponding to input $(d,u) \in A$ satisfies:

$$\left\| T_r(t_1)x \right\|_r > \left\| x_0 \right\|_r + 1 + \frac{1}{\mu(t_1)} a^{-1}\left( R + M\left( t_1, \left\| x_0 \right\|_r + \sup_{t_0 \leq \tau \leq t_1} \left| u(\tau) \right| \right) \right), \quad \forall t \geq t_0 \qquad (2.16)$$

Let $B := \sup_{t_0 \leq \tau \leq t_1} \left\| T_r(\tau)x \right\| < +\infty$. Using (2.8) and (2.15), it follows that (2.14) holds for all $\tilde{x}_0 \in C^1([-r,0]; \Re^n)$ with $\left\| x_0 - \tilde{x}_0 \right\|_r \leq 1$ with $\bar{c} = B + \left\| x_0 \right\|_r + 2 + \max\left\{ \dfrac{a^{-1}(R + M(t,s))}{\mu(t)}; 0 \leq s \leq \left\| x_0 \right\|_r + 1 + \sup_{t_0 \leq \tau \leq t_1} \left| u(\tau) \right|, t_0 \leq t \leq t_1 \right\} + \sup_{t_0 \leq \tau \leq t_1} \left| u(\tau) \right|$, where $\tilde{x}(t)$ denotes the solution of (1.2) with initial condition $T_r(t_0)x = \tilde{x}_0$ corresponding to input $(d,u) \in A$. Let $\varepsilon := \left\| T_r(t_1)x \right\|_r - \left\| x_0 \right\|_r - 1 - \dfrac{1}{\mu(t_1)} a^{-1}\left( R + M\left( t_1, \left\| x_0 \right\|_r + \sup_{t_0 \leq \tau \leq t_1} \left| u(\tau) \right| \right) \right) > 0$. Using (2.15), (2.14), density of $C^1([-r,0]; \Re^n)$ in $C^0([-r,0]; \Re^n)$, continuity of the mapping $x \to g(x) := \left\| x \right\|_r + 1 + \dfrac{1}{\mu(t_1)} a^{-1}\left( R + M\left( t_1, \left\| x \right\|_r + \sup_{t_0 \leq \tau \leq t_1} \left| u(\tau) \right| \right) \right)$, we conclude that there exists $\tilde{x}_0 \in C^1([-r,0]; \Re^n)$ such that:

$$\left\| x_0 - \tilde{x}_0 \right\|_r \leq 1; \quad \left| g(x_0) - g(\tilde{x}_0) \right| \leq \frac{\varepsilon}{2}; \quad \left| \left\| T_r(t_1)x \right\|_r - \left\| T_r(t_1)\tilde{x} \right\|_r \right| \leq \frac{\varepsilon}{2}$$

where $\tilde{x}(t)$ denotes the solution of (1.2) with initial condition $T_r(t_0)x = \tilde{x}_0$ corresponding to input $(d,u) \in A$. Combining (2.15) for $\tilde{x}(t)$ with the above inequalities and the definition of $\varepsilon$ we obtain $\left\| T_r(t_1)x \right\|_r > \left\| T_r(t_1)x \right\|_r$, a contradiction. The proof is complete. ◁

The following definition introduces an important relation between output mappings. The equivalence relation defined next, will be used extensively in the following sections of the present work.

**Definition 2.9:** *Suppose that there exists a continuous mapping* $h : [-r, +\infty) \times \Re^n \to \Re^p$ *with* $h(t,0) = 0$ *for all* $t \geq -r$ *and functions* $a_1, a_2 \in K_\infty$ *such that* $a_1(\left| h(t, x(0)) \right|) \leq \left\| H(t,x) \right\|_Y \leq a_2\left( \sup_{\theta \in [-r,0]} \left| h(t+\theta, x(\theta)) \right| \right)$ *for all* $(t,x) \in \Re^+ \times C^0([-r,0]; \Re^n)$. *Then we say that* $H : \Re^+ \times C^0([-r,0]; \Re^n) \to Y$ *is equivalent to the finite-dimensional mapping* $h$.

For example the identity output mapping $H(t,x) = x \in C^0([-r,0]; \Re^n)$ is equivalent to finite-dimensional mapping $h(t,x) = x \in \Re^n$.

Finally, we end this section by presenting a technical small-gain lemma that will be used in the proofs of our main results. It is a direct corollary of Theorem 1 in [41] and is closely related to Lemma A.1 in [12].

**Lemma 2.10:** *For every $\sigma \in KL$ and $a \in K$ with $a(s) < s$ for all $s > 0$, there exists $\tilde{\sigma} \in KL$ with the following property: if $y : [t_0, t_1) \to \Re^+$, $u : \Re^+ \to \Re^+$ are locally bounded functions and $M \geq 0$ a constant such that the following inequality holds for all $t \in [t_0, t_1)$:*

$$y(t) \leq \inf_{t_0 \leq \xi \leq t} \max \left\{ \sigma(M, t - \xi) \; ; \; a\left( \sup_{\xi \leq \tau \leq t} y(\tau) \right) \; ; \; u(t) \right\} \tag{2.17}$$

*then the following estimate holds for all $t \in [t_0, t_1)$:*

$$y(t) \leq \max \left\{ \tilde{\sigma}(M, t - t_0) \; ; \; \sup_{t_0 \leq \tau \leq t} u(\tau) \right\} \tag{2.18}$$

## 3. Robust Global Asymptotic Output Stability (RGAOS)

In this section we introduce the reader to the notion of non-uniform in time and uniform Robust Global Asymptotic Output Stability (RGAOS) for systems described by RFDEs and we provide different equivalent characterizations for these notions. Notice that the notion of RGAOS is applied to uncertain systems with a robust equilibrium point (vanishing perturbations) and is an "Internal Stability" property.

**Definition 3.1:** *Consider system (1.1) under hypotheses (H1-5). We say that (1.1) is **non-uniformly in time Robustly Globally Asymptotically Output Stable (RGAOS) with disturbances** $d \in M_D$ if (1.1) is RFC and the following properties hold:*

**P1** (1.1) is **Robustly Lagrange Output Stable**, i.e., *for every $\varepsilon > 0$, $T \geq 0$, it holds that*

$$\sup \left\{ \| H(t, \phi(t, t_0, x_0; d)) \|_Y \; ; \; t \in [t_0, +\infty), \; \| x_0 \|_r \leq \varepsilon, \; t_0 \in [0, T], \; d \in M_D \right\} < +\infty$$
**(Robust Lagrange Output Stability)**

**P2** (1.1) is **Robustly Lyapunov Output Stable**, i.e., *for every $\varepsilon > 0$ and $T \geq 0$ there exists a $\delta := \delta(\varepsilon, T) > 0$ such that:*

$$\| x_0 \|_r \leq \delta, \; t_0 \in [0, T] \Rightarrow \| H(t, \phi(t, t_0, x_0; d)) \|_Y \leq \varepsilon, \; \forall t \geq t_0, \; \forall d \in M_D$$
**(Robust Lyapunov Output Stability)**

**P3** (1.1) satisfies the **Robust Output Attractivity Property**, i.e. *for every $\varepsilon > 0$, $T \geq 0$ and $R \geq 0$, there exists a $\tau := \tau(\varepsilon, T, R) \geq 0$, such that:*

$$\| x_0 \|_r \leq R, \; t_0 \in [0, T] \Rightarrow \| H(t, \phi(t, t_0, x_0; d)) \|_Y \leq \varepsilon, \; \forall t \geq t_0 + \tau, \; \forall d \in M_D$$

*Moreover, if there exists a function $a \in K_\infty$ such that $a(\| x \|_r) \leq \| H(t, x) \|_Y$ for all $(t, x) \in \Re^+ \times C^0([-r, 0]; \Re^n)$, then we say that (1.1) is **non-uniformly in time Robustly Globally Asymptotically Stable (RGAS) with disturbances** $d \in M_D$.*

*We say that (1.1) is **non-uniformly in time Robustly Globally Asymptotically Output Stable (RGAOS) with disturbances** $d \in \tilde{M}_D$ if (1.1) is RFC and properties P1-3 above hold with $d \in \tilde{M}_D$ instead of $d \in M_D$.*

The next lemma provides an estimate of the output behavior for non-uniformly in time RGAOS systems. It is an immediate corollary of Lemma 3.4 in [19].

**Lemma 3.2:** *System (1.1) under hypotheses (H1-5) is non-uniformly in time RGAOS with disturbances $d \in M_D$ (or $d \in \tilde{M}_D$) if and only if system (1.1) is RFC and there exist functions $\sigma \in KL$, $\beta \in K^+$ such that the following estimate holds for all $(t_0, x_0) \in \Re^+ \times C^0([-r, 0]; \Re^n)$, $d \in M_D$ (or $d \in \tilde{M}_D$) and $t \geq t_0$:*

$$\|H(t, \phi(t, t_0, x_0; d))\|_Y \leq \sigma\left(\beta(t_0)\|x_0\|_r, t - t_0\right) \tag{3.1}$$

We next provide the definition of Uniform Robust Global Asymptotic Output Stability, in terms of $KL$ functions, which is completely analogous to the finite-dimensional case (see [23,28,39,40]). It is clear that such a definition is equivalent to a $\delta - \varepsilon$ definition (analogous to Definition 3.1).

**Definition 3.3:** *Suppose that (1.1) under hypotheses (H1-5) is RGAOS with disturbances $d \in M_D$ (or $d \in \tilde{M}_D$) and there exist $\sigma \in KL$ such that estimate (3.1) holds for all $(t_0, x_0) \in \Re^+ \times C^0([-r, 0]; \Re^n)$, $d \in M_D$ (or $d \in \tilde{M}_D$) and $t \geq t_0$ with $\beta(t) \equiv 1$. Then we say that (1.1) is **Uniformly Robustly Globally Asymptotically Output Stable (URGAOS) with disturbances** $d \in M_D$ (or $d \in \tilde{M}_D$).*

The following lemma must be compared to Lemma 1.1, page 131 in [8] and Proposition 3.2 in [16]. It shows that for periodic systems RGAOS is equivalent to URGAOS. We say that (1.1) under hypotheses (H1-5) is $T$-periodic, if there exists $T > 0$ such that $f(t + T, x, d) = f(t, x, d)$ and $H(t + T, x) = H(t, x)$ for all $(t, x, d) \in \Re^+ \times C^0([-r, 0]; \Re^n) \times D$. We say that (1.1) under hypotheses (H1-5) is autonomous if $f(t, x, d) = f(0, x, d)$ and $H(t, x) = H(0, x)$ for all $(t, x, d) \in \Re^+ \times C^0([-r, 0]; \Re^n) \times D$.

**Lemma 3.4:** *Suppose that (1.1) under hypotheses (H1-5) is T-periodic. If (1.1) is non-uniformly in time RGAOS with disturbances $d \in M_D$ (or $d \in \tilde{M}_D$), then (1.1) is URGAOS with disturbances $d \in M_D$ (or $d \in \tilde{M}_D$).*

**Proof** The proof is based on the following observation: if (1.1) is $T$-periodic then for all $(t_0, x_0, d) \in \Re^+ \times C^0([-r, 0]; \Re^n) \times M_D$ it holds that $\phi(t, t_0, x_0; d) = \phi(t - kT, t_0 - kT, x_0; P_{kT}d)$ and $H(t, \phi(t, t_0, x_0; d)) = H(t - kT, \phi(t - kT, t_0 - kT, x_0; P_{kT}d))$, where $k := [t_0 / T]$ denotes the integer part of $t_0 / T$ and $(P_{kT}d)(t) = d(t + kT)$ for all $t + kT \geq 0$. Notice that if $d \in M_D$ then $P_{kT}d \in M_D$ and if $d \in \tilde{M}_D$ then $P_{kT}d \in \tilde{M}_D$.

Since (1.1) is non-uniformly in time RGAOS, there exist functions $\sigma \in KL$, $\beta \in K^+$ such that (3.1) holds for all $(t_0, x_0) \in \Re^+ \times C^0([-r, 0]; \Re^n)$, $d \in M_D$ (or $d \in \tilde{M}_D$) and $t \geq t_0$. Consequently, it follows that the following estimate holds for all $(t_0, x_0) \in \Re^+ \times C^0([-r, 0]; \Re^n)$, $d \in M_D$ (or $d \in \tilde{M}_D$) and $t \geq t_0$:

$$\|H(t, \phi(t, t_0, x_0; d))\|_Y \leq \sigma\left(\beta\left(t_0 - \left[\frac{t_0}{T}\right]T\right)\|x_0\|_r, t - t_0\right)$$

Since $0 \leq t_0 - \left[\frac{t_0}{T}\right]T < T$, for all $t_0 \geq 0$, it follows that the following estimate holds for all $(t_0, x_0) \in \Re^+ \times C^0([-r, 0]; \Re^n)$, $d \in M_D$ (or $d \in \tilde{M}_D$) and $t \geq t_0$:

$$\|H(t, \phi(t, t_0, x_0; d))\|_Y \leq \tilde{\sigma}\left(\|x_0\|_r, t - t_0\right)$$

where $\tilde{\sigma}(s, t) := \sigma(Rs, t)$ and $R := \max\{\beta(t); 0 \leq t \leq T\}$. The previous estimate in conjunction with Definition 3.3 implies that (1.1) is URGAOS. The proof is complete. ◁

We are now in a position to state Lyapunov-like characterizations for non-uniform in time RGAOS and URGAOS. The proofs are provided in the Appendix.

**Theorem 3.5** *Consider system (1.1) under hypotheses (H1-5). The following statements are equivalent:*

**(a)** *(1.1) is non-uniformly in time RGAOS with disturbances $d \in M_D$.*

**(b)** *(1.1) is non-uniformly in time RGAOS with disturbances $d \in \tilde{M}_D$.*

**(c)** *(1.1) is RFC and there exist functions $a_1, a_2 \in K_\infty$, $\beta \in K^+$, a positive definite continuous function $\rho : \Re^+ \to \Re^+$ and a mapping $V : \Re^+ \times C^0([-r,0]; \Re^n) \to \Re^+$, which is almost Lipschitz on bounded sets, such that:*

$$a_1\left(\|H(t,x)\|_\gamma\right) \leq V(t,x) \leq a_2\left(\beta(t)\|x\|_r\right), \quad \forall (t,x) \in \Re^+ \times C^0([-r,0]; \Re^n) \tag{3.2}$$

$$V^0(t,x; f(t,x,d)) \leq -\rho(V(t,x)), \quad \forall (t,x,d) \in \Re^+ \times C^0([-r,0]; \Re^n) \times D \tag{3.3}$$

**(d)** *(1.1) is RFC and there exist functions $a_1, a_2 \in K_\infty$, $\beta \in K^+$ and a mapping $V : \Re^+ \times C^0([-r,0]; \Re^n) \to \Re^+$, which is almost Lipschitz on bounded sets, such that inequalities (3.2), (3.3) hold with $\gamma(t) \equiv 1$ and $\rho(s) := s$.*

**(e)** *(1.1) is RFC and there exist a lower semi-continuous mapping $V : \Re^+ \times C^0([-r-\tau,0]; \Re^n) \to \Re^+$, a constant $\tau \geq 0$, functions $a_1, a_2 \in K_\infty$, $\beta, \gamma \in K^+$ with $\int_0^{+\infty} \gamma(t)dt = +\infty$, $\mu \in \mathbf{E}$ (see Notations) and $\rho \in C^0(\Re^+; \Re^+)$ being positive definite, such that the following inequalities hold:*

$$a_1\left(\|H(t,x)\|_\gamma\right) \leq V(t,x) \leq a_2\left(\beta(t)\|x\|_{r+\tau}\right), \quad \forall (t,x) \in \Re^+ \times C^0([-r-\tau,0]; \Re^n) \tag{3.4}$$

$$V^0(t,x; f(t,T_r(0)x,d)) \leq -\gamma(t)\rho(V(t,x)) + \gamma(t)\mu\left(\int_0^t \gamma(s)ds\right), \quad \forall (t,d) \in [\tau,+\infty) \times D, \forall x \in S(t) \tag{3.5}$$

*where the set-valued map $S(t)$ is defined for $t \geq \tau$ by $S(t) := \bigcup_{d \in \tilde{M}_D} S(t,d)$ and the set-valued map $S(t,d)$ is defined for $t \geq \tau$ and $d \in \tilde{M}_D$ by:*

$$S(t,d) := \left\{ x \in C^0([-r-\tau,0]; \Re^n) \,;\, x(\theta) = x(-\tau) + \int_{-\tau}^{\theta} f(t+s, T_r(s)x, d(\tau+s))ds, \forall \theta \in [-\tau,0] \right\} \tag{3.6}$$

*Moreover,*

**i)** *if $H : \Re^+ \times C^0([-r,0]; \Re^n) \to Y$ is equivalent to the finite-dimensional continuous mapping $h : [-r,+\infty) \times \Re^n \to \Re^p$ then inequality (3.2) in statements (c) and (d) can be replaced by the following inequality:*

$$a_1\left(|h(t,x(0))|\right) \leq V(t,x) \leq a_2\left(\beta(t)\|x\|_r\right), \quad \forall (t,x) \in \Re^+ \times C^0([-r,0]; \Re^n) \tag{3.7}$$

**ii)** *if $H : \Re^+ \times C^0([-r,0]; \Re^n) \to Y$ is equivalent to the finite-dimensional continuous mapping $h : [-r,+\infty) \times \Re^n \to \Re^p$ then inequality (3.4) in statement (e) can be replaced by the following inequality:*

$$a_1\left(|h(t,x(0))|\right) \leq V(t,x) \leq a_2\left(\beta(t)\|x\|_{r+\tau}\right), \quad \forall (t,x) \in \Re^+ \times C^0([-r-\tau,0]; \Re^n) \tag{3.8}$$

**iii)** *if there exist functions $a \in K_\infty$, $\mu \in K^+$ and a constant $R \geq 0$ such that $a(\mu(t)|x(0)|) \leq V(t,x) + R$ for all $(t,x) \in \Re^+ \times C^0([-r,0]; \Re^n)$ then the requirement that (1.1) is RFC is not needed in statements (c) and (d) above.*

**Theorem 3.6** *Consider system (1.1) under hypotheses (H1-5). The following statements are equivalent:*

**(a)** *(1.1) is URGAOS with disturbances $d \in M_D$.*

**(b)** *(1.1) is URGAOS with disturbances $d \in \tilde{M}_D$.*

**(c)** *(1.1) is RFC and there exist functions $a_1, a_2 \in K_\infty$, a positive definite continuous function $\rho : \Re^+ \to \Re^+$ and a mapping $V : \Re^+ \times C^0([-r,0];\Re^n) \to \Re^+$, which is almost Lipschitz on bounded sets, such that:*

$$a_1\left(\|H(t,x)\|_Y\right) \leq V(t,x) \leq a_2\left(\|x\|_r\right), \quad \forall (t,x) \in \Re^+ \times C^0([-r,0];\Re^n) \tag{3.9}$$

$$V^0(t,x; f(t,x,d)) \leq -\rho(V(t,x)), \quad \forall (t,x,d) \in \Re^+ \times C^0([-r,0];\Re^n) \times D \tag{3.10}$$

**(d)** *(1.1) is RFC and there exist functions $a_1, a_2 \in K_\infty$ and a mapping $V : \Re^+ \times C^0([-r,0];\Re^n) \to \Re^+$, which is almost Lipschitz on bounded sets, such that inequalities (3.9), (3.10) hold with $\rho(s) := s$. Moreover, if system (1.1) is $T$-periodic, then $V$ is $T$-periodic (i.e. $V(t+T,x) = V(t,x)$ for all $(t,x) \in \Re^+ \times C^0([-r,0];\Re^n)$) and if (1.1) is autonomous then $V$ is independent of $t$.*

**(e)** *(1.1) is RFC and there exist constants $\tau, \beta \geq 0$, a lower semi-continuous mapping $V : \Re^+ \times C^0([-r-\tau,0];\Re^n) \to \Re^+$, functions $a_1, a_2 \in K_\infty$ and $\rho \in C^0(\Re^+;\Re^+)$ being positive definite, such that the following inequalities hold:*

$$a_1\left(\|H(t,x)\|_Y\right) \leq V(t,x) \leq a_2\left(\|x\|_{r+\tau}\right), \quad \forall (t,x) \in \Re^+ \times C^0([-r-\tau,0];\Re^n) \tag{3.11}$$

$$V^0(t,x; f(t,T_r(0)x,d)) \leq \beta V(t,x), \quad \forall (t,x,d) \in \Re^+ \times C^0([-r-\tau,0];\Re^n) \times D \tag{3.12a}$$

$$V^0(t,x; f(t,T_r(0)x,d)) \leq -\rho(V(t,x)), \quad \forall (t,d) \in [\tau,+\infty) \times D, \forall x \in S(t) \tag{3.12b}$$

*where the set-valued map $S(t)$ is defined for $t \geq \tau$ by $S(t) := \bigcup_{d \in \tilde{M}_D} S(t,d)$ and the set-valued map $S(t,d)$ is defined for $t \geq \tau$ and $d \in \tilde{M}_D$ by (3.6).*

*Moreover,*

**i)** *if $H : \Re^+ \times C^0([-r,0];\Re^n) \to Y$ is equivalent to the finite-dimensional continuous $T$-periodic mapping $h : [-r,+\infty) \times \Re^n \to \Re^p$ then inequality (3.9) in statements (c) and (d) can be replaced by the following inequality:*

$$a_1\left(|h(t,x(0))|\right) \leq V(t,x) \leq a_2\left(\|x\|_r\right), \quad \forall (t,x) \in \Re^+ \times C^0([-r,0];\Re^n) \tag{3.13}$$

**ii)** *if $H : \Re^+ \times C^0([-r,0];\Re^n) \to Y$ is equivalent to the finite-dimensional continuous $T$-periodic mapping $h : [-r,+\infty) \times \Re^n \to \Re^p$ then inequality (3.11) in statement (e) can be replaced by the following inequality:*

$$a_1\left(|h(t,x(0))|\right) \leq V(t,x) \leq a_2\left(\|x\|_{r+\tau}\right), \quad \forall (t,x) \in \Re^+ \times C^0([-r-\tau,0];\Re^n) \tag{3.14}$$

**iii)** *if there exist functions $a \in K_\infty$, $\mu \in K^+$ and a constant $R \geq 0$ such that $a(\mu(t)|x(0)|) \leq V(t,x) + R$ for all $(t,x) \in \Re^+ \times C^0([-r,0];\Re^n)$ then the requirement that (1.1) is RFC is not needed in statements (c) and (d) above.*

## 4. Input-to-Output Stability (IOS)

In this section we introduce the reader to the notion of non-uniform in time and uniform Input-to-Output Stability (IOS) for systems described by RFDEs and we provide estimates for the solutions of such systems. Notice that the notion of IOS is an "External Stability" property since it is applied to systems which operate under the effect of external non-vanishing perturbations.

**Definition 4.1:** *We say that (1.2) under hypotheses (S1-7) satisfies the **non-uniform in time Input-to-Output Stability property (IOS) from the input** $u \in M_U$ **with gain** $\gamma \in K$ **and weight** $\delta \in K^+$, if (1.2) is robustly forward complete (RFC) from the input $u \in M_U$ and there exist functions $\sigma \in KL$, $\beta \in K^+$, such that for all $(d,u) \in M_D \times M_U$, $(t_0, x_0) \in \Re^+ \times C^0([-r,0]; \Re^n)$ the solution $x(t)$ of (1.2) with $T_r(t_0)x = x_0$ corresponding to $(d,u) \in M_D \times M_U$ satisfies the following estimate for all $t \geq t_0$:*

$$\|H(t, T_r(t)x)\|_Y \leq \max\left\{ \sigma(\beta(t_0)\|x_0\|_r, t-t_0), \sup_{t_0 \leq \tau \leq t} \gamma(\delta(\tau)|u(\tau)|) \right\} \quad (4.1)$$

*Moreover,*

(i) *if $\beta(t) = \delta(t) \equiv 1$, then we say that (1.2) satisfies the **Uniform Input-to-Output Stability property (UIOS) from the input** $u \in M_U$ **with gain** $\gamma \in K$.*

(ii) *if $\|x\|_r \leq \|H(t,x)\|_Y$ for all $(t,x) \in \Re^+ \times C^0([-r,0]; \Re^n)$, then we say that (1.2) satisfies the **non-uniform in time Input-to-State Stability property (ISS) from the input** $u \in M_U$ **with gain** $\gamma \in K$ **and weight** $\delta \in K^+$.*

(iii) *if $\|x\|_r \leq \|H(t,x)\|_Y$ for all $(t,x) \in \Re^+ \times C^0([-r,0]; \Re^n)$ and $\beta(t) = \delta(t) \equiv 1$, then we say that (1.2) satisfies the **Uniform Input-to-State Stability property (UISS) from the input** $u \in M_U$ **with gain** $\gamma \in K$.*

The following lemma shows that for periodic systems estimate (4.1) leads to a simpler estimate. We say that (1.2) under hypotheses (S1-7) is $T-$periodic, if there exists $T > 0$ such that $f(t+T, x, u, d) = f(t, x, u, d)$ and $H(t+T, x) = H(t, x)$ for all $(t, x, u, d) \in \Re^+ \times C^0([-r,0]; \Re^n) \times U \times D$. We say that (1.2) under hypotheses (S1-7) is autonomous if $f(t, x, d) = f(0, x, d)$ and $H(t, x) = H(0, x)$ for all $(t, x, u, d) \in \Re^+ \times C^0([-r,0]; \Re^n) \times U \times D$.

**Lemma 4.2:** *Suppose that (1.2) is T-periodic. If (1.2) satisfies the non-uniform in time IOS property from the input $u \in M_U$, then there exist functions $\sigma \in KL$, $\delta \in K^+$ and $\gamma \in K$ such that estimate (4.1) holds for all $(t_0, x_0, d, u) \in \Re^+ \times C^0([-r,0]; \Re^n) \times M_D \times M_U$ and $t \geq t_0$ with $\beta(t) \equiv 1$.*

**Proof** The proof is based on the following observation: if (1.2) is $T-$periodic then for all $(t_0, x_0, d, u) \in \Re^+ \times C^0([-r,0]; \Re^n) \times M_D \times M_U$ it holds that $\phi(t, t_0, x_0; u, d) = \phi(t - kT, t_0 - kT, x_0; P_{kT}u, P_{kT}d)$ and $H(t, \phi(t, t_0, x_0; u, d)) = H(t - kT, \phi(t - kT, t_0 - kT, x_0; P_{kT}u, P_{kT}d))$, where $k := [t_0 / T]$ denotes the integer part of $t_0 / T$ and the inputs $P_{kT}u \in M_U$, $P_{kT}d \in M_D$ are defined by $(P_{kT}d)(t) = d(t + kT)$ and $(P_{kT}u)(t) = u(t + kT)$ for all $t + kT \geq 0$.

Since (1.2) satisfies the non-uniform in time IOS property from the input $u \in M_U$, there exist functions $\sigma \in KL$, $\beta, \delta \in K^+$, $\gamma \in K$ such that (4.1) holds for all $(t_0, x_0, d, u) \in \Re^+ \times C^0([-r,0]; \Re^n) \times M_D \times M_U$ and $t \geq t_0$. Consequently, it follows that the following estimate holds for all $(t_0, x_0, d, u) \in \Re^+ \times C^0([-r,0]; \Re^n) \times M_D \times M_U$ and $t \geq t_0$:

$$\|H(t, \phi(t, t_0, x_0; u, d))\|_Y \leq \max\left\{ \sigma(\beta(t_0 - kT)\|x_0\|_r, t-t_0), \sup_{\tau \in [t_0 - kT, t - kT]} \gamma(\delta(\tau)|(P_{kT}u)(\tau)|) \right\}$$

Setting $\tau = s - kT$ and since $0 \leq t_0 - \left[\dfrac{t_0}{T}\right]T < T$, for all $t_0 \geq 0$, we obtain

$$\|H(t,\phi(t,t_0,x_0;u,d))\|_Y \leq \max\left\{\tilde{\sigma}(\|x_0\|_r, t-t_0), \sup_{s\in[t_0,t]} \gamma(\delta(s-kT)|(P_{kT}u)(s-kT)|)\right\} \quad (4.2)$$

where $\tilde{\sigma}(s,t) := \sigma((1+R)s,t)$ and $R := \max\{\beta(t) ; 0 \leq t \leq T\}$. Estimate (4.2) and the identity $(P_{kT}u)(s-kT) = u(s)$ for all $s \geq 0$, imply that the following estimate holds for all $(t_0, x_0, d, u) \in \Re^+ \times C^0([-r,0];\Re^n) \times M_D \times M_U$ and $t \geq t_0$:

$$\|H(t,\phi(t,t_0,x_0;u,d))\|_Y \leq \max\left\{\tilde{\sigma}(\|x_0\|_r, t-t_0), \sup_{s\in[t_0,t]} \gamma(\delta(s-kT)|u(s)|)\right\} \quad (4.3)$$

Setting $\tilde{\delta}(t) := \max\{\delta(s) ; s \in [0,t]\}$, we obtain from (4.3) that estimate (4.1) holds with $\beta(t) \equiv 1$ and $\tilde{\sigma} \in KL$, $\tilde{\delta} \in K^+$ in place of $\sigma \in KL$, $\delta \in K^+$, respectively. The proof is complete. ◁

It should be emphasized that Lemma 4.2 does not guarantee the UIOS property. Particularly, the proof of Lemma 4.2 shows that if (1.2) satisfies the non-uniform in time IOS property with gain $\gamma$ and weight $\delta$ and is $T-periodic$ then estimate (4.1) holds with $\beta(t) \equiv 1$, same gain $\gamma$ but with weight $\tilde{\delta}(t) := \max\{\delta(s) ; s \in [0,t]\}$. Thus the UIOS property can be guaranteed if in addition we assume that $\delta$ is bounded.

We are now in a position to state characterizations for the non-uniform in time IOS property for time-varying uncertain systems. The proof of the following theorem is provided in the Appendix.

**Theorem 4.3:** *The following statements are equivalent for system (1.2) under hypotheses (S1-7):*

**(a)** *System (1.2) is robustly forward complete (RFC) from the input* $u \in M_U$ *and there exist functions* $\sigma \in KL$, $\beta, \phi \in K^+$, $\rho \in K$ *such that that for all* $(d,u) \in M_D \times M_U$, $(t_0, x_0) \in \Re^+ \times C^0([-r,0];\Re^n)$ *the solution* $x(t)$ *of (1.2) with* $T_r(t_0)x = x_0$ *corresponding to* $(d,u) \in M_D \times M_U$ *satisfies the following estimate for all* $t \geq t_0$:

$$\|H(t,T_r(t)x)\|_Y \leq \max\left\{\sigma(\beta(t_0)\|x_0\|_r, t-t_0), \sup_{t_0 \leq \tau \leq t} \sigma(\beta(\tau)\rho(\phi(\tau)|u(\tau)|), t-\tau)\right\} \quad (4.4)$$

**(b)** *System (1.2) satisfies the non-uniform in time IOS property from the input* $u \in M_U$.

**(c)** *There exist a locally Lipschitz function* $\theta \in K_\infty$, *functions* $\phi, \mu \in K^+$ *such that the following system is non-uniformly in time RGAOS with disturbances* $(d',d) \in \tilde{M}_\Delta$:

$$\dot{x}(t) = f\left(t, T_r(t)x, \dfrac{\theta(\|T_r(t)x\|_r)}{\phi(t)}d'(t), d(t)\right) \; ; \; Y(t) = \tilde{H}(t, T_r(t)x) \quad (4.5)$$

where $\Delta := B_U[0,1] \times D$, $\tilde{H}(t,x) := (H(t,x), \mu(t)x) \in Y \times C^0([-r,0];\Re^n)$.

**(d)** *There exist a Lyapunov functional* $V : \Re^+ \times C^0([-r,0];\Re^n) \to \Re^+$, *which is almost Lipschitz on bounded sets, functions* $a_1, a_2, a_3$ *of class* $K_\infty$, $\beta, \delta, \mu$ *of class* $K^+$ *such that:*

$$a_1(\|H(t,x)\|_Y + \mu(t)\|x\|_r) \leq V(t,x) \leq a_2(\beta(t)\|x\|_r), \; \forall (t,x) \in \Re^+ \times C^0([-r,0];\Re^n) \quad (4.6)$$

$$V^0(t,x;f(t,x,u,d)) \leq -V(t,x) + a_3(\delta(t)|u|), \; \forall (t,x,u,d) \in \Re^+ \times C^0([-r,0];\Re^n) \times U \times D \quad (4.7)$$

**(e)** *System (1.2) is RFC from the input* $u \in M_U$ *and there exist a Lyapunov functional* $V: \Re^+ \times C^0([-r,0]; \Re^n) \to \Re^+$, *which is almost Lipschitz on bounded sets, functions* $a_1, a_2, \zeta$ *of class* $K_\infty$, $\beta, \delta$ *of class* $K^+$ *and a continuous positive definite function* $\rho: \Re^+ \to \Re^+$ *such that:*

$$a_1(\|H(t,x)\|_Y) \leq V(t,x) \leq a_2(\beta(t)\|x\|_r), \quad \forall (t,x) \in \Re^+ \times C^0([-r,0]; \Re^n) \tag{4.8}$$

$$V^0(t, x; f(t, x, u, d)) \leq -\rho(V(t,x)),$$
$$\text{for all } (t, x, u, d) \in \Re^+ \times C^0([-r,0]; \Re^n) \times U \times D \text{ with } \zeta(\delta(t)|u|) \leq V(t,x) \tag{4.9}$$

*Moreover,*

**i)** *if* $H: \Re^+ \times C^0([-r,0]; \Re^n) \to Y$ *is equivalent to the finite-dimensional continuous mapping* $h:[-r,+\infty) \times \Re^n \to \Re^p$ *then inequality (4.8) in statement (e) can be replaced by the following inequality:*

$$a_1(|h(t,x(0))|) \leq V(t,x) \leq a_2(\beta(t)\|x\|_r), \quad \forall (t,x) \in \Re^+ \times C^0([-r,0]; \Re^n) \tag{4.10}$$

**ii)** *if there exist functions* $a \in K_\infty$, $\mu \in K^+$ *and a constant* $R \geq 0$ *such that* $a(\mu(t)|x(0)|) \leq \|H(t,x)\|_Y + R$ *for all* $(t,x) \in \Re^+ \times C^0([-r,0]; \Re^n)$ *then the requirement that (1.2) is RFC from the input* $u \in M_U$ *is not needed in statement (a) above.*

**iii)** *if there exist functions* $a \in K_\infty$, $\mu \in K^+$ *and a constant* $R \geq 0$ *such that* $a(\mu(t)|x(0)|) \leq V(t,x) + R$ *for all* $(t,x) \in \Re^+ \times C^0([-r,0]; \Re^n)$ *then the requirement that (1.2) is RFC from the input* $u \in M_U$ *is not needed in statement (e) above.*

In order to obtain characterizations of the uniform IOS property, we need an extra hypothesis for system (1.2).

**(S8)** There exists a constant $R \geq 0$ and a function $a \in K_\infty$ such that the inequality $\|x\|_r \leq a(\|H(t,x)\|_Y) + R$ holds for all $(t,x) \in \Re^+ \times C^0([-r,0]; \Re^n)$.

Hypothesis (S8) holds for the important case of the output map $H(t,x) := d(x(\theta), A); \theta \in [-r,0]$, where $A \subset \Re^n$ is a compact set which contains $0 \in \Re^n$ and $d(x, A)$ denotes the distance of the point $x \in \Re^n$ from the set $A \subset \Re^n$. Notice that it is not required that $A \subset \Re^n$ is positively invariant for (1.2) with $u \equiv 0$.

Hypothesis (S8) allows us to provide characterizations for the UIOS property for periodic uncertain systems. The proof of the following theorem is given in the Appendix.

**Theorem 4.4:** *Suppose that system (1.2) under hypotheses (S1-8) is $T-$ periodic. The following statements are equivalent:*

**(a)** *There exist functions* $\sigma \in KL$, $\rho \in K_\infty$ *such that that for all* $(d,u) \in M_D \times M_U$, $(t_0, x_0) \in \Re^+ \times C^0([-r,0]; \Re^n)$ *the solution* $x(t)$ *of (1.2) with* $T_r(t_0)x = x_0$ *corresponding to* $(d,u) \in M_D \times M_U$, *satisfies the following estimate for all* $t \geq t_0$:

$$\|H(t, T_r(t)x)\|_Y \leq \max\left\{\sigma(\|x_0\|_r, t-t_0), \sup_{t_0 \leq \tau \leq t} \sigma(\rho(|u(\tau)|), t-\tau)\right\} \tag{4.11}$$

**(b)** *System (1.2) satisfies the UIOS property.*

**(c)** *There exists a locally Lipschitz function $\theta \in K_\infty$ such that $0 \in C^0([-r,0]; \Re^n)$ is URGAOS with disturbances $(d',d) \in \tilde{M}_\Delta$ for the system:*

$$\dot{x}(t) = f\left(t, T_r(t)x, \theta\left(\|H(t, T_r(t)x)\|_Y\right)d'(t), d(t)\right) \quad ; \quad Y(t) = H(t, T_r(t)x) \tag{4.12}$$

*where $\Delta := B_U[0,1] \times D$.*

**(d)** *There exist a $T-$ periodic Lyapunov functional $V : \Re^+ \times C^0([-r,0]; \Re^n) \to \Re^+$, which is almost Lipschitz on bounded sets, functions $a_1, a_2, a_3$ of class $K_\infty$ such that:*

$$a_1(\|H(t,x)\|_Y) \leq V(t,x) \leq a_2(\|x\|_r), \quad \forall (t,x) \in \Re^+ \times C^0([-r,0]; \Re^n) \tag{4.13}$$

$$V^0(t,x; f(t,x,u,d)) \leq -V(t,x) + a_3(|u|), \quad \forall (t,x,u,d) \in \Re^+ \times C^0([-r,0]; \Re^n) \times U \times D \tag{4.14}$$

**(e)** *There exist a Lyapunov functional $V : \Re^+ \times C^0([-r,0]; \Re^n) \to \Re^+$, which is almost Lipschitz on bounded sets, functions $a_1, a_2, \zeta$ of class $K_\infty$ and a continuous positive definite function $\rho : \Re^+ \to \Re^+$ such that:*

$$a_1(\|H(t,x)\|_Y) \leq V(t,x) \leq a_2(\|x\|_r), \quad \forall (t,x) \in \Re^+ \times C^0([-r,0]; \Re^n) \tag{4.15}$$

$$V^0(t,x; f(t,x,u,d)) \leq -\rho(V(t,x)),$$
$$\text{for all } (t,x,u,d) \in \Re^+ \times C^0([-r,0]; \Re^n) \times U \times D \text{ with } \zeta(|u|) \leq V(t,x) \tag{4.16}$$

*Finally, if $H : \Re^+ \times C^0([-r,0]; \Re^n) \to Y$ is equivalent to the finite-dimensional continuous $T-$ periodic mapping $h : [-r, +\infty) \times \Re^n \to \Re^p$ then inequalities (4.13), (4.15) in statements (d) and (e), respectively, can be replaced by the following inequality:*

$$a_1(|h(t, x(0))|) \leq V(t,x) \leq a_2(\|x\|_r), \quad \forall (t,x) \in \Re^+ \times C^0([-r,0]; \Re^n) \tag{4.17}$$

**Remark 4.5:** A Statement like (e) of Theorem 4.4 was extensively used as a tool of proving the uniform ISS property for autonomous time-delay systems in [33].

The following theorem provides sufficient Lyapunov-like conditions for the non-uniform in time and uniform IOS property. The proof of implications (e) $\Rightarrow$ (a) of Theorem 4.3 and (e) $\Rightarrow$ (a) of Theorem 4.4 are based on the result of Theorem 4.6, which gives quantitative estimates of the solutions of (1.2) under hypotheses (S1-7). The gain functions and the weights of the IOS property can be determined explicitly in terms of the functions involved in the assumptions of Theorem 4.6.

**Theorem 4.6:** *Consider system (1.2) under hypotheses (S1-7) and suppose that there exist a Lyapunov functional $V : \Re^+ \times C^0([-r,0]; \Re^n) \to \Re^+$, which is almost Lipschitz on bounded sets, functions $a, \zeta$ of class $K_\infty$, $\beta, \delta$ of class $K^+$ and a continuous positive definite function $\rho : \Re^+ \to \Re^+$ such that:*

$$V(t,x) \leq a(\beta(t)\|x\|_r), \quad \forall (t,x) \in \Re^+ \times C^0([-r,0]; \Re^n) \tag{4.18}$$

$$V^0(t,x; f(t,x,u,d)) \leq -\rho(V(t,x)),$$
$$\text{for all } (t,x,u,d) \in \Re^+ \times C^0([-r,0]; \Re^n) \times U \times D \text{ with } \zeta(\delta(t)|u|) \leq V(t,x) \tag{4.19}$$

*Moreover, suppose that one of the following holds:*
**a)** *system (1.2) is RFC from the input $u \in M_U$*

**b)** *there exist functions $p \in K_\infty$, $\mu \in K^+$ and a constant $R \geq 0$ such that $p(\mu(t)|x(0)|) \leq V(t,x) + R$ for all $(t,x) \in \Re^+ \times C^0([-r,0]; \Re^n)$*

Then system (1.2) is RFC from the input $u \in M_U$ and there exist a function $\sigma \in KL$ with $\sigma(s,0) = s$ for all $s \geq 0$, such that that for all $(d,u) \in M_D \times M_U$, $(t_0, x_0) \in \Re^+ \times C^0([-r,0]; \Re^n)$ the solution $x(t)$ of (1.2) with $T_r(t_0)x = x_0$ corresponding to $(d,u) \in M_D \times M_U$, satisfies the following estimate for all $t \geq t_0$:

$$V(t, T_r(t)x) \leq \max\left\{ \sigma\big(a(\beta(t_0)\|x_0\|_r), t - t_0 \big), \sup_{t_0 \leq \tau \leq t} \sigma\big(\zeta(\delta(\tau)|u(\tau)|), t - \tau \big) \right\} \quad (4.20)$$

*Finally,*

(i) *if there exist a function $a_1$ of class $K_\infty$ such that $a_1(\|H(t,x)\|_Y) \leq V(t,x)$ for all $(t,x) \in \Re^+ \times C^0([-r,0]; \Re^n)$, then system (1.2) satisfies the non-uniform in time IOS property from the input $u \in M_U$ with gain $\gamma(s) := a_1^{-1}(\zeta(s))$ and weight $\delta$.*

(ii) *if $H : \Re^+ \times C^0([-r,0]; \Re^n) \to Y$ is equivalent to the finite-dimensional continuous mapping $h : [-r, +\infty) \times \Re^n \to \Re^p$ and there exist functions $a_1, a_2$ of class $K_\infty$ such that $a_1(|h(t, x(0))|) \leq V(t,x)$, $\|H(t,x)\|_Y \leq a_2\left(\sup_{\theta \in [-r,0]} |h(t+\theta, x(\theta))|\right)$ for all $(t,x) \in \Re^+ \times C^0([-r,0]; \Re^n)$, then system (1.2) satisfies the non-uniform in time IOS property from the input $u \in M_U$ with gain $\gamma(s) := a_2(a_1^{-1}(\zeta(s)))$ and weight $\delta$.*

The proof of Theorem 4.6 is based on the following comparison lemma. Its proof is given in the Appendix.

**Lemma 4.7:** *For each positive definite continuous function $\rho : \Re^+ \to \Re^+$ there exists a function $\sigma$ of class $KL$, with $\sigma(s,0) = s$ for all $s \geq 0$ with the following property: if $y : [t_0, t_1] \to \Re^+$ is an absolutely continuous function, $u : \Re^+ \to \Re^+$ is a locally bounded mapping and $I \subset [t_0, t_1]$ a set of Lebesgue measure zero such that $\dot{y}(t)$ is defined on $[t_0, t_1] \setminus I$ and such that the following implication holds for all $t \in [t_0, t_1] \setminus I$:*

$$y(t) \geq u(t) \quad \Rightarrow \quad \dot{y}(t) \leq -\rho(y(t)) \quad (4.21)$$

*then the following estimate holds for all $t \in [t_0, t_1]$:*

$$y(t) \leq \max\left\{ \sigma(y(t_0), t - t_0), \sup_{t_0 \leq s \leq t} \sigma(u(s), t - s) \right\} \quad (4.22)$$

The following example presents an autonomous time-delay system, which satisfies the non-uniform in time IOS property and does not satisfy the UIOS property. The analysis is performed with the help of Theorem 4.6.

**Example 4.8:** Consider the following autonomous time-delay system:

$$\begin{aligned}
\dot{x}_1(t) &= d(t)x_1(t) \\
\dot{x}_2(t) &= -x_2(t) + x_1(t-r)u(t) \\
Y(t) &= x_2(t) \\
(x_1(t), x_2(t))' &\in \Re^2, d(t) \in D := [-1,1], u(t) \in U := \Re, Y(t) \in \Re
\end{aligned} \quad (4.23)$$

Consider the functional:

$$V(t, x_1, x_2) := \exp(-8t)x_1^4(0) + \exp(-4t)x_1^2(0) + \frac{1}{2}x_2^2(0) + \frac{1}{4}\exp(-8t)\int_{-r}^{0} x_1^4(s)ds \qquad (4.24)$$

First notice that the functional $V : \Re^+ \times C^0([-r,0];\Re^2) \to \Re^+$ defined by (4.24) is almost Lipschitz on bounded sets. Moreover, inequality (4.18) is satisfied for this functional with $a(s) := \left(1 + \frac{r}{4}\right)s^4 + s^2$ and $\beta(t) \equiv 1$. We next estimate an upper bound for the Dini derivative of the functional $V : \Re^+ \times C^0([-r,0];\Re^2) \to \Re^+$ along the solutions of system (4.23). We have for all $(t, x, u, d) \in \Re^+ \times C^0([-r,0];\Re^2) \times \Re \times [-1,1]$:

$$V^0(t, x_1, x_2; (dx_1(0), -x_2(0) + x_1(-r)u)) = -8\exp(-8t)x_1^4(0) + 4d\exp(-8t)x_1^4(0) - 4\exp(-4t)x_1^2(0) + 2d\exp(-4t)x_1^2(0)$$
$$- x_2^2(0) + x_2(0)x_1(-r)u - 2\exp(-8t)\int_{-r}^{0} x_1^4(s)ds + \frac{1}{4}\exp(-8t)x_1^4(0) - \frac{1}{4}\exp(-8t)x_1^4(-r)$$

Using the inequalities $|d| \leq 1$, $x_2(0)x_1(-r)u \leq \frac{1}{2}x_2^2(0) + \frac{1}{2}x_1^2(-r)u^2$, $\frac{1}{2}x_1^2(-r)u^2 \leq \frac{1}{4}\exp(-8t)x_1^4(-r) + \frac{1}{4}\exp(8t)u^4$, we are in a position to estimate for all $(t, x, u, d) \in \Re^+ \times C^0([-r,0];\Re^2) \times \Re \times [-1,1]$:

$$V^0(t, x_1, x_2; (dx_1(0), -x_2(0) + x_1(-r)u)) \leq$$
$$-3\exp(-8t)x_1^4(0) - 2\exp(-4t)x_1^2(0) - \frac{1}{2}x_2^2(0) + \frac{1}{4}\exp(8t)u^4 - 2\exp(-8t)\int_{-r}^{0} x_1^4(s)ds$$

Finally, using the above inequality and definition (4.24) we obtain for all $(t, x, u, d) \in \Re^+ \times C^0([-r,0];\Re^2) \times \Re \times [-1,1]$:

$$V^0(t, x_1, x_2; (dx_1(0), -x_2(0) + x_1(-r)u)) \leq -V(t, x_1, x_2) + \frac{1}{4}\exp(8t)u^4 \qquad (4.25)$$

Inequality (4.25) guarantees that (4.19) holds with $\rho(s) := \frac{1}{2}s$, $\zeta(s) := \frac{1}{2}s^4$ and $\delta(t) := \exp(2t)$. Definition (4.24) guarantees that there exist functions $p \in K_\infty$, $\mu \in K^+$ and a constant $R \geq 0$ such that $p(\mu(t)|x(0)|) \leq V(t,x) + R$ for all $(t,x) \in \Re^+ \times C^0([-r,0];\Re^2)$ (for instance, $p(s) := \frac{1}{2}s^2$, $\mu(t) := \exp(-2t)$ and $R := 0$). Finally, notice that the inequality $a_1(\|H(t,x)\|_Y) \leq V(t,x)$ holds for all $(t,x) \in \Re^+ \times C^0([-r,0];\Re^2)$ with $a_1(s) := \frac{1}{2}s^2$. It follows from Theorem 4.6 that system (4.23) satisfies the non-uniform in time IOS property from the input $u \in M_U$ with gain $\gamma(s) := a_1^{-1}(\zeta(s)) = s^2$ and weight $\delta(t) := \exp(2t)$.

It should be emphasized that system (4.23) does not satisfy the UIOS property from the input $u \in M_U$. This can be shown by considering the solution of (4.23) corresponding to inputs $d(t) \equiv 1$ and $u(t) \equiv 1$. It can be shown that for $x_1(0) \neq 0$ the output of (4.23) is not bounded and satisfies $\lim_{t \to +\infty} |Y(t)| = +\infty$. Consequently, bounded inputs can produce unbounded output responses, which contradicts the requirements of the UIOS property from the input $u \in M_U$. ◁

# 5. Razumikhin Functions

Let $V:[-r,+\infty)\times\Re^n \to \Re$ be a locally Lipschitz mapping and let $(t,x,v)\in \Re^+ \times \Re^n \times \Re^n$. We define

$$D^+V(t,x;v) := \limsup_{h\to 0^+} \frac{V(t+h,x+hv)-V(t,x)}{h} \tag{5.1}$$

The following proposition extends the classical Razumikhin theorem to systems with disturbances as well as to the case of output asymptotic stability. Its proof is based on Lemma 2.10 (small-gain lemma) and on Lemma 4.7 of the previous section.

**Proposition 5.1:** *Consider system (1.1) under hypotheses (H1-5) and suppose that $H:\Re^+\times C^0([-r,0];\Re^n)\to Y$ is equivalent to the finite-dimensional mapping $h:[-r,+\infty)\times\Re^n \to \Re^p$. Moreover, suppose that there exist a locally Lipschitz function $V:[-r,+\infty)\times\Re^n \to \Re^+$, functions $a_1$, $a_2$, $a$ of class $K_\infty$, with $a(s)<s$ for all $s>0$, $\beta$ of class $K^+$ and a positive definite function $\rho$ such that:*

$$a_1(|h(t-r,x)|) \le V(t-r,x) \le a_2(\beta(t)|x|), \quad \forall (t,x)\in \Re^+ \times \Re^n \tag{5.2}$$

$$D^+V(t,x(0);f(t,x,d)) \le -\rho(V(t,x(0))),$$
$$\text{for all } (t,x,d)\in \Re^+ \times C^0([-r,0];\Re^n)\times D \text{ with } a\left(\sup_{\theta\in[-r,0]} V(t+\theta,x(\theta))\right) \le V(t,x(0)) \tag{5.3}$$

*Finally, suppose that one of the following holds:*

*(i) system (1.1) is RFC*

*(ii) there exist functions $\zeta\in K_\infty$, $\mu\in K^+$ and a constant $R\ge 0$ such that $\zeta(\mu(t)|x|)\le V(t-r,x)+R$ for all $(t,x)\in \Re^+\times \Re^n$*

*Then system (1.1) is non-uniformly in time RGAOS. Moreover, if $\beta$ is bounded then system (1.1) is URGAOS.*

**Proof:** Consider a solution $x(t)$ of (1.1) under hypotheses (H1-5) corresponding to arbitrary $d\in M_D$ with initial condition $T_r(t_0)x = x_0 \in C^0([-r,0];\Re^n)$. It follows from (5.3) and Lemma 4.7 that there exists a continuous function $\sigma$ of class $KL$, with $\sigma(s,0)=s$ for all $s\ge 0$ such that:

$$V(t,x(t)) \le \max\left\{\sigma(V(t_0,x(t_0)),t-t_0),\sup_{t_0\le s\le t}\sigma\left(a\left(\sup_{\theta\in[-r,0]}V(s+\theta,x(s+\theta))\right),t-s\right)\right\}, \quad \forall t\in[t_0,t_{max}) \tag{5.4}$$

An immediate consequence of estimate (5.4) and the fact that $\sigma(s,0)=s$ for all $s\ge 0$ is the following estimate for all $t\in[t_0,t_{max})$:

$$\sup_{\theta\in[-r,0]} V(t+\theta,x(t+\theta)) \le \max\left\{\overline{\sigma}\left(\sup_{\theta\in[-r,0]}V(t_0+\theta,x(t_0+\theta)),t-t_0\right); a\left(\sup_{t_0\le s\le t}\sup_{\theta\in[-r,0]}V(s+\theta,x(s+\theta))\right)\right\} \tag{5.5}$$

where $\overline{\sigma}(s,t):=s$ for $t\in[0,r]$ and $\overline{\sigma}(s,t):=\sigma(s,t-r)$ for $t>r$. Using the fact that $a(s)<s$ for all $s>0$ and estimate (5.5) it may be shown that:

$$\sup_{\theta\in[-r,0]}V(t+\theta,x(t+\theta)) \le \sup_{\theta\in[-r,0]}V(t_0+\theta,x(t_0+\theta)), \quad \forall t\in[t_0,t_{max}) \tag{5.6}$$

In case that system (1.1) is RFC, we have $t_{max} = +\infty$. In case that there exist functions $\zeta \in K_\infty$, $\mu \in K^+$ and a constant $R \geq 0$ such that $\zeta(\mu(t)|x|) \leq V(t-r, x) + R$ for all $(t, x) \in \Re^+ \times \Re^n$, inequality (5.6) in conjunction with (5.2) implies:

$$\|T_r(t)x\|_r \leq \frac{1}{\min_{0 \leq \tau \leq t+r} \mu(t)} \zeta^{-1}\left(R + a_2\left(\max_{0 \leq \tau \leq t_0+r} \beta(\tau)\|x_0\|_r\right)\right), \quad \forall t \in [t_0, t_{max}) \tag{5.7}$$

Estimate (5.7) implies that system (1.1) is RFC. Therefore we conclude that in any case system (1.1) is RFC and that estimates (5.5) and (5.6) holds for all $t \geq t_0$. Combining (5.5) with (5.6) we obtain for all $t \geq t_0$:

$$\sup_{\theta \in [-r,0]} V(t+\theta, x(t+\theta)) \leq \inf_{t_0 \leq \xi \leq t} \max\left\{\bar{\sigma}\left(\sup_{\theta \in [-r,0]} V(t_0+\theta, x(t_0+\theta)), t-\xi\right); a\left(\sup_{\xi \leq s \leq t} \sup_{\theta \in [-r,0]} V(s+\theta, x(s+\theta))\right)\right\} \tag{5.8}$$

Lemma 2.10 in conjunction with inequality (5.8) implies the existence of $\tilde{\sigma} \in KL$ such that:

$$\sup_{\theta \in [-r,0]} V(t+\theta, x(t+\theta)) \leq \tilde{\sigma}\left(\sup_{\theta \in [-r,0]} V(t_0+\theta, x(t_0+\theta)), t-t_0\right), \quad \forall t \geq t_0 \tag{5.9}$$

Since $H : \Re^+ \times C^0([-r,0]; \Re^n) \to Y$ is equivalent to the finite-dimensional mapping $h$, there exist function $a_3 \in K_\infty$ such that $\|H(t,x)\|_Y \leq a_3\left(\sup_{\theta \in [-r,0]} |h(t+\theta, x(\theta))|\right)$ for all $(t,x) \in \Re^+ \times C^0([-r,0]; \Re^n)$. Using the previous inequality in conjunction with (5.9) and (5.2) we obtain:

$$\|H(t, T_r(t)x)\|_Y \leq a_3\left(a_1^{-1}\left(\tilde{\sigma}\left(a_2\left(\max_{0 \leq \tau \leq t_0+r} \beta(\tau)\|x_0\|_r\right), t-t_0\right)\right)\right), \quad \forall t \geq t_0 \tag{5.10}$$

Estimate (5.10) implies that system (1.1) is non-uniformly in time RGAOS. Moreover, if $\beta$ is bounded then estimate (5.10) implies that system (1.1) is URGAOS. The proof is complete. ◁

The following example illustrates the use of Proposition 5.1 to a robust stabilization problem for a linear time-varying control system with distributed delay.

**Example 5.2:** Consider the following time-varying control system:

$$\begin{aligned}
\dot{x}_1(t) &= d(t)\exp(t)\int_{t-r}^{t} x_1(s)ds + x_2(t) \\
\dot{x}_2(t) &= u(t) \\
Y(t) &= T_r(t)x \in C^0([-r,0]; \Re^2) \\
x(t) &= (x_1(t), x_2(t))' \in \Re^2, d(t) \in D := [-1,1], u(t) \in \Re
\end{aligned} \tag{5.11}$$

where $r \geq 0$, and the quadratic function:

$$V(t,x) := \frac{33}{4}\exp(2t)x_1^2 + \frac{1}{2}(x_2 + 4\exp(t)x_1)^2 \tag{5.12}$$

First notice that the output $H(t,x) = x \in C^0([-r,0]; \Re^2)$ is equivalent to the finite-dimensional mapping $[-r, +\infty) \times \Re^n \ni (t,x) \to h(t,x) = x \in \Re^2$. Completing the squares, it may be shown that $V$ satisfies (5.2) with $h(x) := (x_1, x_2)' \in \Re^2$, $a_1(s) := \frac{1}{4}s^2$, $\beta(t) := \exp(t)$ and $a_2(s) := 30s^2$. Next we evaluate the derivative of $V$ along the solutions of system (5.11). We obtain for all $(t, x, u, d) \in \Re^+ \times C^0([-r,0]; \Re^2) \times \Re \times [-1,1]$:

$$D^+V\left(t, x(0); \left(d\exp(t)\int_{-r}^0 x_1(s)ds + x_2(0), u\right)'\right) = \frac{33}{2}\exp(2t)x_1^2(0) + \frac{33}{2}d\exp(3t)x_1(0)\int_{-r}^0 x_1(s)ds - 66\exp(3t)x_1^2(0)$$

$$+ (x_2(0) + 4\exp(t)x_1(0))\left(u + 4\exp(t)x_1(0) + \frac{33}{2}\exp(2t)x_1(0) + 4d\exp(2t)\int_{-r}^0 x_1(s)ds + 4\exp(t)x_2(0)\right)$$

Let $a(s) := \frac{1}{2}s$. Thus the inequality $a\left(\sup_{\theta \in [-r,0]} V(t+\theta, x(\theta))\right) \leq V(t, x(0))$ gives

$\left|\int_{-r}^0 x_1(s)ds\right| \leq r\exp(r)\sqrt{2}|x_1(0)| + \frac{2}{\sqrt{33}}r\exp(r)\exp(-t)|x_2(0) + 4\exp(t)x_1(0)|$. Using the previous inequality as well as the fact that $|d| \leq 1$, we obtain for all $(t, x, u, d) \in \Re^+ \times C^0([-r,0]; \Re^2) \times \Re \times [-1,1]$ with $a\left(\sup_{\theta \in [-r,0]} V(t+\theta, x(\theta))\right) \leq V(t, x(0))$:

$$D^+V\left(t, x(0); \left(d\exp(t)\int_{-r}^0 x_1(s)ds + x_2(0), u\right)'\right) \leq -\left(66 - \frac{33}{2} - \frac{33}{\sqrt{2}}K\right)\exp(3t)x_1^2(0)$$

$$+ \left(\sqrt{33} + 4\sqrt{2}\right)K\exp(2t)|x_1(0)||x_2(0) + 4\exp(t)x_1(0)| + \frac{8K}{\sqrt{33}}\exp(t)(x_2(0) + 4\exp(t)x_1(0))^2$$

$$+ (x_2(0) + 4\exp(t)x_1(0))\left(u + 4\exp(t)x_1(0) + \frac{33}{2}\exp(2t)x_1(0) + 4\exp(t)x_2(0)\right)$$

where $K := r\exp(r)$. Consequently, if

$$\frac{3\sqrt{2}}{2} > r\exp(r) \tag{5.13}$$

there exists $\varepsilon > 0$ such that $c := \frac{99}{2} - \left(\frac{33}{\sqrt{2}} + \frac{\sqrt{33} + 4\sqrt{2}}{2\varepsilon}\right)K > 0$. It follows from (5.13) that the linear time-varying feedback law:

$$u(t) = -4\exp(t)x_1(t) - \frac{33}{2}\exp(2t)x_1(t) - 4\exp(t)x_2(t) - L\exp(t)(x_2(t) + 4\exp(t)x_1(t)) \tag{5.14}$$

where $L \geq \left(\frac{8}{\sqrt{33}} + \frac{\varepsilon(\sqrt{33} + 4\sqrt{2})}{2}\right)K + c$, guarantees the following inequality for all $(t, x, d) \in \Re^+ \times C^0([-r,0]; \Re^2) \times [-1,1]$ with $a\left(\sup_{\theta \in [-r,0]} V(t+\theta, x(\theta))\right) \leq V(t, x(0))$:

$$D^+V\left(t, x(0); \left(d\exp(t)\int_{-r}^0 x_1(s)ds + x_2(0), u\right)'\right) \leq -\frac{4c}{33}\exp(t)V(t, x(0))$$

Hence, it follows from Proposition 5.1 that the closed-loop system (5.11), (5.14) under hypothesis (5.13) is non-uniformly in time RGAS. ◁

The following proposition provides conditions in terms of Razumikhin functions for the non-uniform in time and uniform IOS property.

**Proposition 5.3:** *Consider system (1.2) under hypotheses (S1-7) and suppose that $H : \mathfrak{R}^+ \times C^0([-r,0];\mathfrak{R}^n) \to Y$ is equivalent to the finite-dimensional mapping $h:[-r,+\infty) \times \mathfrak{R}^n \to \mathfrak{R}^p$. Moreover, suppose that there exist a locally Lipschitz function $V:[-r,+\infty) \times \mathfrak{R}^n \to \mathfrak{R}^+$, functions $a_1, a_2, a, \zeta$ of class $K_\infty$ with $a(s) < s$ for all $s > 0$, $\beta, \delta$ of class $K^+$ and a positive definite function $\rho$ such that:*

$$a_1(|h(t-r,x)|) \leq V(t-r,x) \leq a_2(\beta(t)|x|), \quad \forall (t,x) \in \mathfrak{R}^+ \times \mathfrak{R}^n \quad (5.15)$$

$$D^+V(t,x(0); f(t,x,u,d)) \leq -\rho(V(t,x(0))),$$

*for all $(t,x,u,d) \in \mathfrak{R}^+ \times C^0([-r,0];\mathfrak{R}^n) \times U \times D$ with $\max\left\{\zeta(\delta(t)|u|), a\left(\sup_{\theta \in [-r,0]} V(t+\theta, x(\theta))\right)\right\} \leq V(t,x(0))$* (5.16)

*Finally, suppose that one of the following holds:*

(i) *system (1.1) is RFC from the input $u \in M_U$*

(ii) *there exist functions $p \in K_\infty$, $\mu \in K^+$ and a constant $R \geq 0$ such that $p(\mu(t)|x|) \leq V(t-r,x) + R$ for all $(t,x) \in \mathfrak{R}^+ \times \mathfrak{R}^n$*

*Let $a_3 \in K_\infty$ the function with the property $\|H(t,x)\|_Y \leq a_3\left(\sup_{\theta \in [-r,0]} |h(t+\theta, x(\theta))|\right)$ for all $(t,x) \in \mathfrak{R}^+ \times C^0([-r,0];\mathfrak{R}^n)$. Then system (1.2) satisfies the non-uniform in time IOS property with gain $\gamma(s) := a_3(a_1^{-1}(\zeta(s)))$ and weight $\delta$. Moreover, if $\beta, \delta$ are bounded then system (1.2) satisfies the UIOS property from the input $u \in M_U$.*

**Proof:** Consider a solution $x(t)$ of (1.2) under hypotheses (S1-7) corresponding to arbitrary $(u,d) \in M_U \times M_D$ with initial condition $T_r(t_0)x = x_0 \in C^0([-r,0];\mathfrak{R}^n)$. It follows from (5.16) and Lemma 4.7 that there exists a continuous function $\sigma$ of class $KL$, with $\sigma(s,0) = s$ for all $s \geq 0$ such that for all $t \in [t_0, t_{max})$ we have:

$$V(t,x(t)) \leq \max\left\{\sigma(V(t_0,x(t_0)), t-t_0); \sup_{t_0 \leq s \leq t} \sigma\left(a\left(\sup_{\theta \in [-r,0]} V(s+\theta, x(s+\theta))\right), t-s\right); \sup_{t_0 \leq s \leq t} \sigma(\zeta(\delta(s)|u(s)|), t-s)\right\} \quad (5.17)$$

An immediate consequence of estimate (5.17) and the fact that $\sigma(s,0) = s$ for all $s \geq 0$ is the following estimate for all $t \in [t_0, t_{max})$:

$$\sup_{\theta \in [-r,0]} V(t+\theta, x(t+\theta)) \leq \max\left\{\bar{\sigma}\left(\sup_{\theta \in [-r,0]} V(t_0+\theta, x(t_0+\theta)), t-t_0\right); a\left(\sup_{t_0 \leq s \leq t} \sup_{\theta \in [-r,0]} V(s+\theta, x(s+\theta))\right); \sup_{t_0 \leq s \leq t} \zeta(\delta(s)|u(s)|)\right\} \quad (5.18)$$

where $\bar{\sigma}(s,t) := s$ for $t \in [0,r]$ and $\bar{\sigma}(s,t) := \sigma(s, t-r)$ for $t > r$. Using the fact that $a(s) < s$ for all $s > 0$ and estimate (5.18) it may be shown that:

$$\sup_{\theta \in [-r,0]} V(t+\theta, x(t+\theta)) \leq \max\left\{\sup_{\theta \in [-r,0]} V(t_0+\theta, x(t_0+\theta)); \sup_{t_0 \leq s \leq t} \zeta(\delta(s)|u(s)|)\right\}, \forall t \in [t_0, t_{max}) \quad (5.19)$$

In case that system (1.2) is RFC from the input $u \in M_U$, we have $t_{\max} = +\infty$. In case that there exist functions $p \in K_\infty$, $\mu \in K^+$ and a constant $R \geq 0$ such that $p(\mu(t)|x|) \leq V(t-r, x) + R$ for all $(t, x) \in \Re^+ \times \Re^n$, inequality (5.19) in conjunction with (5.15) implies:

$$\|T_r(t)x\|_r \leq \frac{1}{\min_{0 \leq \tau \leq t+r} \mu(t)} p^{-1}\left(R + a_2\left(\max_{0 \leq \tau \leq t_0+r} \beta(\tau)\|x_0\|_r\right) + \sup_{t_0 \leq s \leq t} \zeta(\delta(s)|u(s)|)\right), \forall t \in [t_0, t_{\max}) \quad (5.20)$$

Estimate (5.20) implies that system (1.2) is RFC from the input $u \in M_U$. Therefore we conclude that in any case system (1.2) is RFC from the input $u \in M_U$ and that estimates (5.18), (5.19) hold for all $t \geq t_0$. Combining (5.18) with (5.19) we obtain for all $t \geq t_0$:

$$\sup_{\theta \in [-r,0]} V(t+\theta, x(t+\theta)) \leq$$
$$\inf_{t_0 \leq \xi \leq t} \max\left\{\bar{\sigma}\left(\sup_{\theta \in [-r,0]} V(t_0+\theta, x(t_0+\theta)), t-\xi\right) ; a\left(\sup_{\xi \leq s \leq t} \sup_{\theta \in [-r,0]} V(s+\theta, x(s+\theta))\right) ; \sup_{t_0 \leq s \leq t} \zeta(\delta(s)|u(s)|)\right\} \quad (5.21)$$

Lemma 2.10 in conjunction with inequality (5.21) implies the existence of $\tilde{\sigma} \in KL$ such that:

$$\sup_{\theta \in [-r,0]} V(t+\theta, x(t+\theta)) \leq \max\left\{\tilde{\sigma}\left(\sup_{\theta \in [-r,0]} V(t_0+\theta, x(t_0+\theta)), t-t_0\right) ; \sup_{t_0 \leq s \leq t} \zeta(\delta(s)|u(s)|)\right\}, \forall t \geq t_0 \quad (5.22)$$

Since $H : \Re^+ \times C^0([-r,0]; \Re^n) \to Y$ is equivalent to the finite-dimensional mapping $h$, there exists a function $a_3 \in K_\infty$ such that $\|H(t,x)\|_Y \leq a_3\left(\sup_{\theta \in [-r,0]} |h(t+\theta, x(\theta))|\right)$ for all $(t,x) \in \Re^+ \times C^0([-r,0]; \Re^n)$. Using the previous inequality in conjunction with (5.22) and (5.15) we obtain:

$$\|H(t, T_r(t)x)\|_Y \leq \max\left\{a_3\left(a_1^{-1}\left(\tilde{\sigma}\left(a_2\left(\max_{0 \leq \tau \leq t_0+r} \beta(\tau)\|x_0\|_r\right), t-t_0\right)\right)\right) ; \sup_{t_0 \leq s \leq t} a_3\left(a_1^{-1}(\zeta(\delta(s)|u(s)|))\right)\right\}, \forall t \geq t_0$$
(5.23)

Estimate (5.23) implies that system (1.2) satisfies the non-uniform in time IOS property from the input $u \in M_U$ with gain $\gamma(s) := a_3(a_1^{-1}(\zeta(s)))$ and weight $\delta$. Moreover, if $\beta, \delta$ are bounded then estimate (5.23) implies that system (1.2) satisfies the UIOS property from the input $u \in M_U$. The proof is complete. ◁

The following example illustrates the application of Proposition 5.3 to an autonomous time-delay system.

**Example 5.4:** Consider the following autonomous time-delay system:

$$\begin{aligned}\dot{x}(t) &= d(t)x(t-r) - x^3(t) + u(t) \\ Y(t) &= H(T_r(t)x) \\ x(t) &\in \Re, d(t) \in D := [-R, R], u(t) \in U := \Re, Y(t) \in C^0([-r,0]; \Re)\end{aligned} \quad (5.24)$$

where $R > 0$, $H(x) := h(x(\theta)); \theta \in [-r,0]$, $h(x) := x\left(1 - 2\sqrt{R}|x|^{-1}\right)$ for $|x| > 2\sqrt{R}$ and $h(x) := 0$ for $|x| \leq 2\sqrt{R}$. Notice that $H : C^0([-r,0]; \Re) \to Y := C^0([-r,0]; \Re)$ is equivalent to the finite-dimensional mapping $h : \Re \to \Re$. Consider the locally Lipschitz function:

$$V(x) := \max\{0; x^2 - 4R\} \quad (5.25)$$

which satisfies (5.15) with $a_1(s) = a_2(s) := s^2$ and $\beta(t) \equiv 1$. Notice that for all $|d| \leq R$ and $x \in C^0([-r,0]; \Re)$ with $|x(0)| > 2\sqrt{R}$, we have:

$$D^+V(x(0); dx(-r) - x^3(0) + u) = 2dx(0)x(-r) - 2x^4(0) + 2x(0)u \leq 2R|x(0)||x(-r)| - 2x^4(0) + 2|x(0)||u|$$

Using the Young inequality $2|x(0)||u| \leq x^4(0) + 3|u|^{\frac{4}{3}}$ and completing the squares, we obtain for all $|d| \leq R$ and $x \in C^0([-r,0]; \Re)$ with $|x(0)| > 2\sqrt{R}$:

$$D^+V(x(0); dx(-r) - x^3(0) + u) \leq 2R\, x^2(0) + \frac{R}{2} x^2(-r) - x^4(0) + 3|u|^{\frac{4}{3}} \quad (5.26)$$

Let $a(s) := \frac{1}{4}s$, $\rho(s) := 2Rs$ and $\zeta(s) := \frac{3}{2R} s^{\frac{4}{3}}$. Notice that if $a\left(\sup_{\theta \in [-r,0]} V(x(\theta))\right) \leq V(x(0))$ and $|x(0)| > 2R$ then $x^2(-r) \leq 4x^2(0)$. Consequently, previous definitions, definition (5.25) and inequality (5.26) implies that the following inequality holds for all $|d| \leq R$ and $x \in C^0([-r,0]; \Re)$ with $\max\left\{\zeta(|u|), a\left(\sup_{\theta \in [-r,0]} V(x(\theta))\right)\right\} \leq V(x(0))$, $|x(0)| > 2\sqrt{R}$:

$$D^+V(x(0); dx(-r) - x^3(0) + u) \leq -\rho(V(x(0))) \quad (5.27)$$

Notice that (5.27) holds also for all $|d| \leq R$ and $x \in C^0([-r,0]; \Re)$ with $\max\left\{\zeta(|u|), a\left(\sup_{\theta \in [-r,0]} V(x(\theta))\right)\right\} \leq V(x(0))$, $|x(0)| \leq 2\sqrt{R}$. Thus $V$ satisfies (5.16) with $a(s) := \frac{1}{4}s$, $\rho(s) := 2Rs$, $\zeta(s) := \frac{3}{2R} s^{\frac{4}{3}}$ and $\delta(t) \equiv 1$. Finally, notice that there exist functions $p \in K_\infty$, $\mu \in K^+$ and a constant $K \geq 0$ such that $p(\mu(t)|x|) \leq V(t-r, x) + K$ for all $(t,x) \in \Re^+ \times \Re^n$ (for instance, $p(s) := s^2$, $\mu(t) \equiv 1$ and $K = 4R$). It follows from Proposition 5.3 that system (5.24) satisfies the UIOS property from the input $u \in M_U$ with gain $\gamma(s) := \sqrt{\frac{3}{2R}} s^{\frac{2}{3}}$ and weight $\delta(t) \equiv 1$. ◁

## 6. Conclusions

In this work Lyapunov-like characterizations of non-uniform in time and uniform Robust Global Asymptotic Output Stability (RGAOS) for uncertain time-varying systems described by Retarded Functional Differential Equations (RFDEs) are developed. Moreover, the notions of uniform and non-uniform in time Input-to-State Stability (ISS) and Input-to-Output Stability (IOS) are introduced. Necessary and sufficient conditions in terms of Lyapunov functionals and Razumikhin functions are provided for these notions. The framework of the present work allows outputs with no delays, outputs with discrete or distributed delays and functional outputs with memory.

The robust stability notions and properties proposed in the present work are parallel to those recently developed for dynamical systems described by finite-dimensional ordinary differential equations. Just as the popularity gained by the notions of uniform and non-uniform in time RGAOS, ISS and IOS in the context of deterministic systems, it is our firm belief that the stability results of this paper will play an important role in mathematical systems and control theory for important classes of systems described by RFDEs. We expect to report on our future findings along this direction elsewhere.

**Acknowledgments:** The work of I. Karafyllis has been supported by the State Scholarships Foundation (I.K.Y.). The work of P. Pepe has been supported by the Italian MIUR Project PRIN 2005. The work of Z. P. Jiang has been supported in part by the U.S. National Science Foundation under grants ECS-0093176, OISE-0408925 and DMS-0504462.

## Appendix

**Proof of Theorem 3.5:** Implications (a) $\Rightarrow$ (b), (d) $\Rightarrow$ (c), (c) $\Rightarrow$ (e) are obvious. Thus we are left with the proof of implications (b) $\Rightarrow$ (d), (c) $\Rightarrow$ (a) and (e) $\Rightarrow$ (b).

Proof of (b) $\Rightarrow$ (d): The proof of this implication is based on the methodology presented in [1] for finite-dimensional systems as well as the methodologies followed in [22,28].

Since (1.1) is non-uniformly in time RGAOS with disturbances $d \in \tilde{M}_D$, there exist functions $\sigma \in KL$, $\beta \in K^+$ such that estimate (3.1) holds for all $(t_0, x_0, d) \in \Re^+ \times C^0([-r,0];\Re^n) \times \tilde{M}_D$ and $t \geq t_0$. Moreover, by recalling Proposition 7 in [38] there exist functions $\tilde{a}_1$, $\tilde{a}_2$ of class $K_\infty$, such that the $KL$ function $\sigma(s,t)$ is dominated by $\tilde{a}_1^{-1}(\exp(-2t)\tilde{a}_2(s))$. Thus, by taking into account estimate (3.1), we have:

$$\tilde{a}_1\left(\left\|H(t,\phi(t,t_0,x_0;d))\right\|_Y\right) \leq \exp(-2(t-t_0))\tilde{a}_2\left(\beta(t_0)\|x_0\|_r\right), \forall t \geq t_0 \geq 0, x_0 \in C^0([-r,0];\Re^n), d \in \tilde{M}_D \quad (A1)$$

Without loss of generality we may assume that $\tilde{a}_1 \in K_\infty$ is globally Lipschitz on $\Re^+$ with unit Lipschitz constant, namely, $|\tilde{a}_1(s_1) - \tilde{a}_1(s_2)| \leq |s_1 - s_2|$ for all $s_1, s_2 \geq 0$. To see this notice that we can always replace $\tilde{a}_1 \in K_\infty$ by the function $\bar{a}_1(s) := \inf\left\{\min\left\{\frac{1}{2}y, a(y)\right\} + |y-s| ; y \geq 0\right\}$, which is of class $K_\infty$, globally Lipschitz on $\Re^+$ with unit Lipschitz constant and satisfies $\bar{a}_1(s) \leq \tilde{a}_1(s)$. Moreover, without loss of generality we may assume that $\beta \in K^+$ is non-decreasing.

Since (1.1) is Robustly Forward Complete (RFC), by virtue of Lemma 3.5 in [19], there exist functions $\mu \in K^+$, $a \in K_\infty$ and a constant $M \geq 0$, such that for every $(t_0, x_0, d) \in \Re^+ \times C^0([-r,0];\Re^n) \times M_D$, we have:

$$\|\phi(t,t_0,x_0;d)\|_r \leq \mu(t)\, a(\|x_0\|_r + M), \quad \forall t \geq t_0 \tag{A2}$$

Moreover, without loss of generality we may assume that $\mu \in K^+$ is non-decreasing. Making use of (2.4) and (A2), we obtain the following elementary property for the solution of (1.1):

$$\|H(t,\phi(t,t_0,x;d)) - H(t,\phi(t,t_0,y;d))\|_Y \leq B(t,\|x\|_r + \|y\|_r) \exp\left(\tilde{L}(t,\|x\|_r + \|y\|_r)(t-t_0)\right)\|x-y\|_r \tag{A3}$$

for all $t \geq t_0$ and $(t_0,x,y,d) \in \Re^+ \times C^0([-r,0];\Re^n) \times C^0([-r,0];\Re^n) \times \tilde{M}_D$

where

$$\tilde{L}(t,s) := L(t, 2\mu(t)a(s+M)) \ ; \ B(t,s) := L_H(t, 2\mu(t)a(s+M))$$

and $L(\cdot)$, $L_H(\cdot)$ are the functions involved in (2.1) and (2.2). Furthermore, under hypotheses (H1-4), Lemma 3.2 in [19] implies the existence of functions $\zeta \in K_\infty$ and $\gamma \in K^+$ such that:

$$|f(t,x,d)| \leq \zeta(\gamma(t)\|x\|_r), \quad \forall (t,x,d) \in \Re^+ \times C^0([-r,0];\Re^n) \times D$$

Without loss of generality, we may assume that $\gamma \in K^+$ is non-decreasing. Since $x(t) = x(0) + \int_{t_0}^{t} f(\tau, T_r(\tau)x, d(\tau))d\tau$, using the previous inequality in conjunction with (A2) we obtain:

$$|x(t) - x(0)| \leq (t-t_0) G_1(t,\|x\|_r)$$
$$G_1(t,s) := \zeta(\gamma(t)\mu(t)a(s+M))$$

and consequently

$$\|\phi(t,t_0,x;d) - x\|_r \leq (t-t_0) G_1(t,\|x\|_r) + G_2(x, t-t_0) \tag{A4}$$

for all $t \geq t_0$ and $(t_0,x,d) \in \Re^+ \times C^0([-r,0];\Re^n) \times \tilde{M}_D$

where the functional

$$G_2(x,h) := \sup\{|x(0) - x(\theta)|; \theta \in [-\min(h,r),0]\} + \begin{cases} 0 & \text{if } h \geq r \\ \sup\{|x(\theta+h) - x(\theta)|; \theta \in [-r,-h]\} & \text{if } 0 \leq h < r \end{cases}$$

is defined for all $(x,h) \in C^0([-r,0];\Re^n) \times \Re^+$. Notice that $\lim_{h \to 0^+} G_2(x,h) = 0$ for all $x \in C^0([-r,0];\Re^n)$ and consequently for every $\varepsilon > 0$, $R \geq 0$, $x \in C^0([-r,0];\Re^n)$, there exists $T(\varepsilon,R,x) > 0$ such that:

$$t_0 \leq t \leq t_0 + T(\varepsilon,R,x) \Rightarrow \|\phi(t,t_0,x;d) - x\|_r \leq \varepsilon, \text{ for all } (t_0,x,d) \in [0,R] \times C^0([-r,0];\Re^n) \times \tilde{M}_D \tag{A5}$$

We define for all $q \in N$:

$$U_q(t,x) := \sup\left\{\max\left\{0, \tilde{a}_1(\|H(\tau,\phi(\tau,t,x;d))\|_Y) - q^{-1}\right\} \exp((\tau-t)) : \tau \geq t, d \in \tilde{M}_D\right\} \tag{A6}$$

Clearly, estimate (A1) and definition (A6) imply that:

$$\max\left\{0, \tilde{a}_1(\|H(t,x)\|_Y) - q^{-1}\right\} \leq U_q(t,x) \leq \tilde{a}_2(\beta(t)\|x\|_r), \quad \forall (t,x,q) \in \Re^+ \times C^0([-r,0];\Re^n) \times N \tag{A7}$$

Moreover, by definition (A6) we obtain for all $(h,t,x,d,q) \in \Re^+ \times \Re^+ \times C^0([-r,0];\Re^n) \times \tilde{M}_D \times N$:

$$U_q(t+h, \phi(t+h,t,x;d)) \leq \exp(-h) U_q(t,x) \tag{A8}$$

By virtue of estimate (A1) it follows that for every $(q,R) \in N \times \Re^+$, $\tau \geq t + \tilde{T}(R,q)$, $(t,d) \in [0,R] \times \tilde{M}_D$, and $x \in C^0([-r,0];\Re^n)$ with $\|x\|_r \leq R$, it holds: $\tilde{a}_1(\|H(\tau,\phi(\tau,t,x;d))\|_Y) \leq \exp(-2(\tau-t))\tilde{a}_2(\beta(t)\|x\|_r) \leq q^{-1}$, where

$$\tilde{T}(R,q) := \max\left\{0, \frac{1}{2}\log(1+q\tilde{a}_2(\beta(R)R))\right\} \tag{A9}$$

Thus, by virtue of definitions (A6), (A9), we conclude that:

$$U_q(t,x) = \sup\{\max\{0, \tilde{a}_1(\|H(\tau,\phi(\tau,t,x;d))\|_Y) - q^{-1}\}\exp((\tau-t)) : t \leq \tau \leq t+\xi, d \in \tilde{M}_D\}, \forall \xi \geq \tilde{T}(\max\{t,\|x\|_r\},q) \tag{A10}$$

It follows by taking into account (A10) that for all $t \in [0, R]$, and $(x, y) \in C^0([-r,0]; \Re^n) \times C^0([-r,0]; \Re^n)$ with $\|x\|_r \leq R$, $\|y\|_r \leq R$, it holds:

$$\begin{aligned}
&|U_q(t,y) - U_q(t,x)| = \\
&|\sup\{\max\{0, \tilde{a}_1(\|H(\tau,\phi(\tau,t,y;d))\|_Y) - q^{-1}\}\exp((\tau-t)) : t \leq \tau \leq t+\tilde{T}(R,q), d \in \tilde{M}_D\} \\
&- \sup\{\max\{0, \tilde{a}_1(\|H(\tau,\phi(\tau,t,x;d))\|_Y) - q^{-1}\}\exp((\tau-t)) : t \leq \tau \leq t+\tilde{T}(R,q), d \in \tilde{M}_D\}| \\
&\leq \sup\{\exp((\tau-t))|\tilde{a}_1(\|H(\tau,\phi(\tau,t,y;d))\|_Y) - \tilde{a}_1(\|H(\tau,\phi(\tau,t,y;d))\|_Y)| : t \leq \tau \leq t+\tilde{T}(R,q), d \in \tilde{M}_D\} \\
&\leq \sup\{\exp((\tau-t))\|H(\tau,\phi(\tau,t,y;d)) - H(\tau,\phi(\tau,t,x;d))\|_Y : t \leq \tau \leq t+\tilde{T}(R,q), d \in \tilde{M}_D\}
\end{aligned} \tag{A11}$$

Notice that in the above inequalities we have used the facts that the functions $\max\{0, s - q^{-1}\}$ and $\tilde{a}_1(s)$ are globally Lipschitz on $\Re^+$ with unit Lipschitz constant. From (A3) and (A11) we deduce for all $t \in [0, R]$, and $(x, y) \in C^0([-r,0]; \Re^n) \times C^0([-r,0]; \Re^n)$ with $\|x\|_r \leq R$, $\|y\|_r \leq R$:

$$|U_q(t,y) - U_q(t,x)| \leq G_3(R,q)\|y - x\|_r \tag{A12}$$

where

$$G_3(R,q) := B(R + \tilde{T}(R,q), 2R)\exp(\tilde{T}(R,q)(1 + \tilde{L}(R + \tilde{T}(R,q), 2R))) \tag{A13}$$

Next, we establish continuity with respect to $t$ on $\Re^+ \times C^0([-r,0]; \Re^n)$. Let $R \geq 0$, $q \in N$ arbitrary, $t_1, t_2 \in [0, R]$ with $t_1 \leq t_2$, and $x \in C^0([-r,0]; \Re^n)$ with $\|x\|_r \leq R$. Clearly, we have for all $d \in \tilde{M}_D$:

$$|U_q(t_1,x) - U_q(t_2,x)| \leq (1 - \exp(-(t_2-t_1)))U_q(t_1,x) + |\exp(-(t_2-t_1))U_q(t_1,x) - U_q(t_2,\phi(t_2,t_1,x;d))|$$
$$+ |U_q(t_2,\phi(t_2,t_1,x;d)) - U_q(t_2,x)|$$

By virtue of (A4), (A5), (A8), (A12) and the previous inequality we obtain for all $t_1, t_2 \in [0, R]$ with $t_1 \leq t_2 \leq t_1 + T(1, R, x)$ (where $T(\varepsilon, R, x) > 0$ is involved in (A5)) and $d \in \tilde{M}_D$:

$$\begin{aligned}
|U_q(t_1,x) - U_q(t_2,x)| &\leq (t_2-t_1)U_q(t_1,x) + \exp(-(t_2-t_1))U_q(t_1,x) - U_q(t_2,\phi(t_2,t_1,x;d)) \\
&+ G_3(R+1,q)[G_2(x,t_2-t_1) + (t_2-t_1)G_1(R,R)]
\end{aligned} \tag{A14}$$

Definition (A6) implies that for every $\varepsilon > 0$, there exists $d_\varepsilon \in \tilde{M}_D$ with the following property:

$$U_q(t_1,x) - \varepsilon \leq \sup\{\max\{0, \tilde{a}_1(\|H(\tau,\phi(\tau,t_1,x;d_\varepsilon))\|_Y) - q^{-1}\}\exp((\tau-t_1)); \tau \geq t_1\} \leq U_q(t_1,x) \tag{A15}$$

Thus using definition (A6) we obtain:

$$\begin{aligned}
&\exp(-(t_2-t_1))U_q(t_1,x) - U_q(t_2,\phi(t_2,t_1,x;d_\varepsilon)) \\
&\leq \max\{A_q(t_1,t_2,x), B_q(t_1,t_2,x)\} - B_q(t_1,t_2,x) + \varepsilon\exp(-(t_2-t_1))
\end{aligned} \tag{A16}$$

where

$$A_q(t_1,t_2,x) := \sup\{\max\{0, \tilde{a}_1(\|H(\tau,\phi(\tau,t_1,x;d_\varepsilon))\|_Y) - q^{-1}\}\exp((\tau-t_2)); t_2 \geq \tau \geq t_1\}$$
$$B_q(t_1,t_2,x) := \sup\{\max\{0, \tilde{a}_1(\|H(\tau,\phi(\tau,t_1,x;d_\varepsilon))\|_Y) - q^{-1}\}\exp((\tau-t_2)); \tau \geq t_2\}$$
(A17)

Since the functions $\max\{0, s-q^{-1}\}$ and $\tilde{a}_1(s)$ are globally Lipschitz on $\Re^+$ with unit Lipschitz constant, we obtain:

$$\begin{aligned}&A_q(t_1,t_2,x) - B_q(t_1,t_2,x) \\ &\leq \sup\{\max\{0, \tilde{a}_1(\|H(\tau,\phi(\tau,t_1,x;d_\varepsilon))\|_Y) - q^{-1}\}\exp((\tau-t_2)); t_2 \geq \tau \geq t_1\} - \max\{0, \tilde{a}_1(\|H(t_2,\phi(t_2,t_1,x;d_\varepsilon))\|_Y) - q^{-1}\} \\ &\leq \sup\{\max\{0, \tilde{a}_1(\|H(\tau,\phi(\tau,t_1,x;d_\varepsilon))\|_Y) - q^{-1}\}; t_2 \geq \tau \geq t_1\} - \max\{0, \tilde{a}_1(\|H(t_2,\phi(t_2,t_1,x;d_\varepsilon))\|_Y) - q^{-1}\} \\ &\leq \sup\{|\tilde{a}_1(\|H(\tau,\phi(\tau,t_1,x;d_\varepsilon))\|_Y) - \tilde{a}_1(\|H(t_2,\phi(t_2,t_1,x;d_\varepsilon))\|_Y)|; t_2 \geq \tau \geq t_1\} \\ &\leq \sup\{\|H(\tau,\phi(\tau,t_1,x;d_\varepsilon)) - H(t_2,\phi(t_2,t_1,x;d_\varepsilon))\|_Y; t_2 \geq \tau \geq t_1\}\end{aligned}$$
(A18)

Notice that by virtue of (2.2), (A4) and (A5), we obtain for all $\tau \in [t_1,t_2]$ with $t_1 \leq t_2 \leq t_1 + T(1,R,x)$, $t_1,t_2 \in [0,R]$:

$$\begin{aligned}\|H(\tau,\phi(\tau,t_1,x;d_\varepsilon)) - H(t_2,\phi(t_2,t_1,x;d_\varepsilon))\|_Y &\leq \|H(\tau,\phi(\tau,t_1,x;d_\varepsilon)) - H(t_1,x)\|_Y + \|H(t_2,\phi(t_2,t_1,x;d_\varepsilon)) - H(t_1,x)\|_Y \\ &\leq 2(t_2-t_1)L_H(R,2R+2)(1+G_1(R,R)) + 2L_H(R,2R+2)\sup\{G_2(x,h); h \in [0,t_2-t_1]\}\end{aligned}$$
(A19)

Distinguishing the cases $A_q(t_1,t_2,x) \geq B_q(t_1,t_2,x)$ and $A_q(t_1,t_2,x) \leq B_q(t_1,t_2,x)$ it follows from (A16), (A17), (A18) and (A19) that:

$$\begin{aligned}&\exp(-(t_2-t_1))U_q(t_1,x) - U_q(t_2,\phi(t_2,t_1,x;d_\varepsilon)) \leq \\ &\leq 2(t_2-t_1)L_H(R,2R+2)(1+G_1(R,R)) + 2L_H(R,2R+2)\sup\{G_2(x,h); h \in [0,t_2-t_1]\} + \varepsilon\end{aligned}$$

Combining the previous inequality with (A14) and the right hand side of (A7), we obtain:

$$\begin{aligned}&|U_q(t_1,x) - U_q(t_2,x)| \\ &\leq (t_2-t_1)(\tilde{a}_2(\beta(R)R) + G_1(R,R)G_3(R+1,q) + 2L_H(R,2R+2)(1+G_1(R,R))) \\ &+ (2L_H(R,2R+2) + G_3(R+1,q))\sup\{G_2(x,h); h \in [0,t_2-t_1]\} + \varepsilon\end{aligned}$$
(A20)

Since (A20) holds for all $\varepsilon > 0$, $R \geq 0$, $q \in N$, $x \in C^0([-r,0];\Re^n)$ with $\|x\|_r \leq R$ and $t_1,t_2 \in [0,R]$ with $t_1 \leq t_2 \leq t_1 + T(1,R,x)$, it follows that:

$$|U_q(t_1,x) - U_q(t_2,x)| \leq G_4(R,q)\left[|t_2-t_1| + \sup\{G_2(x,h); h \in [0,|t_2-t_1|]\}\right]$$

for all $R \geq 0$, $q \in N$, $x \in C^0([-r,0];\Re^n)$ with $\|x\|_r \leq R$ and $t_1,t_2 \in [0,R]$ with $|t_2-t_1| \leq T(1,R,x)$ (A21)

where $G_4(R,q) := \tilde{a}_2(\beta(R)R) + (1+G_1(R,R))G_3(R+1,q) + 2L_H(R,2R+2)(1+G_1(R,R))$. Finally, we define:

$$V(t,x) := \sum_{q=1}^{\infty} \frac{2^{-q}U_q(t,x)}{1+G_3(q,q)+G_4(q,q)}$$
(A22)

Inequality (A7) in conjunction with definition (A22) implies (3.2) with $a_2 = \tilde{a}_2$ and $a_1(s) := \sum_{q=1}^{\infty} \frac{2^{-q}\max\{0, a_1(s)-q^{-1}\}}{1+G_3(q,q)+G_4(q,q)}$, which is a function of class $K_\infty$. Moreover, by virtue of definition (A22) and inequality (A8) we obtain for all $(h,t,x,d) \in \Re^+ \times \Re^+ \times C^0([-r,0];\Re^n) \times \tilde{M}_D$:

$$V(t+h,\phi(t+h,t,x;d)) \leq \exp(-h)V(t,x)$$
(A23)

Next define

$$M(R) := 1 + \sum_{q=1}^{[R]+1} \frac{2^{-q} G_3(R,q)}{1 + G_3(q,q) + G_4(q,q)} \tag{A24}$$

which is a positive non-decreasing function. Using (A12) and definition (A22) as well as the fact $G_3(R,q) \leq G_3(q,q)$ for $q > R$, we conclude that property (P1) of Definition 2.5 holds. Let $d \in D$ and define $\tilde{d}(t) \equiv d$. Definition (2.10) and inequality (A23) implies that for all $(t,x) \in \Re^+ \times C^0([-r,0]; \Re^n)$ we have:

$$V^0(t,x; f(t,x,d)) := \limsup_{\substack{h \to 0^+ \\ y \to 0, y \in C^0([-r,0]; \Re^n)}} \frac{V(t+h, E_h(x; f(t,x,d)) + hy) - V(t,x)}{h}$$

$$\leq \limsup_{h \to 0^+} \frac{V(t+h, \phi(t+h,t,x;\tilde{d})) - V(t,x)}{h} + \limsup_{\substack{h \to 0^+ \\ y \to 0, y \in C^0([-r,0]; \Re^n)}} \frac{V(t+h, E_h(x; f(t,x,d)) + hy) - V(t+h, \phi(t+h,t,x;\tilde{d}))}{h}$$

$$\leq -V(t,x) + \limsup_{\substack{h \to 0^+ \\ y \to 0, y \in C^0([-r,0]; \Re^n)}} \frac{V(t+h, E_h(x; f(t,x,d)) + hy) - V(t+h, \phi(t+h,t,x;\tilde{d}))}{h}$$

Let $R \geq \max\{t, \|x\|_r\}$. Definition (2.9) and property (A5) implies that $t + h \leq R + 1$, $\|\phi(t+h,t,x;\tilde{d})\|_r \leq R+1$, $\|E_h(x; f(t,x,d)) + hy\|_r \leq R+1$ for $h$ and $\|y\|_r$ sufficiently small. Using property (P1) of Definition 2.5 and the previous inequalities we obtain:

$$V^0(t,x; f(t,x,d)) \leq -V(t,x) + M(R+1) \limsup_{h \to 0^+} \frac{\|E_h(x; f(t,x,d)) - \phi(t+h,t,x;\tilde{d})\|_r}{h}$$

We set $\phi(t+h,t,x;\tilde{d}) = x(t+h+\theta)$; $\theta \in [-r,0]$. Notice that $\phi(t+h,t,x;\tilde{d}) - E_h(x; f(t,x,d)) = h y_h$, where

$$y_h := \begin{cases} \dfrac{\theta + h}{h} \left( \dfrac{x(t+\theta+h) - x(t)}{\theta+h} - f(t,x,d) \right) & \text{for } -h < \theta \leq 0 \\ 0 & \text{for } -r \leq \theta \leq -h \end{cases}$$

with $\|y_h\|_r \leq \sup\left\{ \left| \dfrac{x(t+s) - x(t)}{s} - f(t,x,d) \right|; 0 < s \leq h \right\}$. Since $\lim_{h \to 0^+} \dfrac{x(t+h) - x(t)}{h} = f(t,x,d)$, we obtain that $y_h \to 0$ as $h \to 0^+$. Hence, we obtain $\limsup_{h \to 0^+} \dfrac{\|E_h(x; f(t,x,d)) - \phi(t+h,t,x;\tilde{d})\|_r}{h} = 0$ and consequently (3.3) holds with $\gamma(t) \equiv 1$ and $\rho(s) := s$.

Finally, we establish continuity of $V$ with respect to $t$ on $\Re^+ \times C^0([-r,0]; \Re^n)$ and property (P2) of Definition 2.5. Notice that by virtue of (A21) and the fact $G_4(R,q) \leq G_4(q,q)$ for $q > R$, we obtain:

$$|V(t_1, x) - V(t_2, x)| \leq P(R) \left[ |t_2 - t_1| + \sup\{G_2(x,h); h \in [0, |t_2 - t_1|]\} \right]$$
for all $R \geq 0$, $x \in C^0([-r,0]; \Re^n)$ with $\|x\|_r \leq R$ and $t_1, t_2 \in [0,R]$ with $|t_2 - t_1| \leq T(1,R,x)$ \quad (A25)

where

$$P(R) := 1 + \sum_{q=1}^{[R]+1} \frac{2^{-q} G_4(R,q)}{1 + G_3(q,q) + G_4(q,q)}$$

is a positive non-decreasing function. Inequality (A25) in conjunction with the fact that $\lim_{h \to 0^+} G_2(x,h) = 0$ for all $x \in C^0([-r,0]; \Re^n)$, establishes continuity of $V$ with respect to $t$ on $\Re^+ \times C^0([-r,0]; \Re^n)$. Moreover, for every absolutely continuous function $x:[-r,0] \to \Re^n$ with $\|x\|_r \leq R$ and essentially bounded derivative, it holds that $\sup\{G_2(x,h) ; h \in [0, |t_2 - t_1|]\} \leq |t_2 - t_1| \sup_{-r \leq \tau \leq 0} |\dot{x}(\tau)|$ for $|t_2 - t_1| \leq r$. It follows from (A4), (A5) and the previous inequality that $T(1,R,x) \geq \dfrac{r}{(1+r)(1+G_1(R,R) + \sup_{-r \leq \tau \leq 0} |\dot{x}(\tau)|)}$. Property (P2) of Definition 2.5 with $G(R) := \dfrac{1+r}{r}(1 + G_1(R,R) + R)$ is a direct consequence of (A25) and the two previous inequalities.

Proof of (c) $\Rightarrow$ (a):

Case 1: (1.1) is RFC

Consider a solution $x(t)$ of (1.1) under hypotheses (H1-5) corresponding to arbitrary $d \in M_D$ with initial condition $T_r(t_0)x = x_0 \in C^1([-r,0]; \Re^n)$. By virtue of Lemma 2.6, for every $T \in (t_0, +\infty)$, the mapping $[t_0, T] \ni t \to V(t, T_r(t)x)$ is absolutely continuous. It follows from (3.3) and Lemma 2.4 that $\dfrac{d}{dt}(V(t, T_r(t)x)) \leq -\rho(V(t, T_r(t)x))$ a.e. on $[t_0, +\infty)$. The previous differential inequality in conjunction with Lemma 4.4 in [28] shows that there exists $\sigma \in KL$ such that

$$V(t, T_r(t)x) \leq \sigma(V(t_0, x_0), t - t_0) \text{ for all } t \geq t_0 \tag{A26}$$

It follows from Lemma 2.7 that the solution $x(t)$ of (1.1) under hypotheses (H1-5) corresponding to arbitrary $d \in M_D$ with arbitrary initial condition $T_r(t_0)x = x_0 \in C^0([-r,0]; \Re^n)$ satisfies (A26) for all $t \geq t_0$. Next, we distinguish the following cases:

1) If (3.2) holds, then (3.1) is a direct consequence of (A26) and (3.2).
2) If (3.7) holds, then (A26) implies the following estimate:

$$|h(t, x(t))| \leq a_1^{-1}\left(\sigma\left(a_2\left(\beta(t_0)\|x_0\|_r\right), t - t_0\right)\right), \forall t \geq t_0$$

Since $h:[-r, +\infty) \times \Re^n \to \Re^p$ is continuous with $h(t, 0) = 0$ for all $t \geq -r$, it follows from Lemma 3.2 in [19] that there exist functions $\zeta \in K_\infty$ and $\gamma \in K^+$ such that:

$$|h(t - r, x)| \leq \zeta(\gamma(t)|x|), \forall (t, x) \in \Re^+ \times \Re^n$$

Combining the two previous inequalities we obtain:

$$\sup_{\theta \in [-r,0]} |h(t + \theta, x(t + \theta))| \leq \omega\left(\phi(t_0)\|x_0\|_r, t - t_0\right), \forall t \geq t_0$$

where $\omega(s,t) := \max\{\zeta(s), a_1^{-1}(\sigma(a_2(s), 0))\}$ for $t \in [0, r)$, $\omega(s,t) := \max\{\exp(r-t)\zeta(s), a_1^{-1}(\sigma(a_2(s), t - r))\}$ for $t \geq r$ and $\phi(t) := \beta(t) + \max_{0 \leq \tau \leq t + r} \gamma(\tau)$. The above estimate, in conjunction with the fact that $H: \Re^+ \times C^0([-r,0]; \Re^n) \to Y$ is equivalent to the finite-dimensional mapping $h$ shows that (1.1) is non-uniformly RGAOS with disturbances $d \in M_D$.

**Case 2:** There exist functions $a \in K_\infty$, $\mu \in K^+$ and a constant $R \geq 0$ such that

$$a(\mu(t)|x(0)|) \leq V(t,x) + R \text{ for all } (t,x) \in \Re^+ \times C^0([-r,0];\Re^n) \quad (A27)$$

Consider a solution of (1.1) under hypotheses (H1-5) corresponding to arbitrary $d \in M_D$ with initial condition $T_r(t_0)x = x_0 \in C^1([-r,0];\Re^n)$. By virtue of Lemma 2.6, for every $T \in (t_0, t_{\max})$, the mapping $[t_0,T] \ni t \to V(t,T_r(t)x)$ is absolutely continuous. It follows from (3.3) and Lemma 2.4 that for every $T \in (t_0, t_{\max})$ it holds that $\frac{d}{dt}(V(t,T_r(t)x)) \leq -\rho(V(t,T_r(t)x))$ a.e. on $[t_0,T]$. The previous differential inequality in conjunction with Lemma 4.4 in [28] shows that there exists $\sigma \in KL$ such that

$$V(t,T_r(t)x) \leq \sigma(V(t_0,x_0),t-t_0) \text{ for all } t \in [t_0,T] \quad (A28)$$

Combining (A27), (A28) and (3.2) we obtain:

$$|x(t)| \leq \frac{1}{\mu(t)} a^{-1}(\sigma(V(t_0,x_0),t-t_0) + R) \text{ for all } t \in [t_0,T] \quad (A29)$$

Estimate (A29) shows that $t_{\max} = +\infty$ and consequently estimates (A28), (A29) hold for all $t \geq t_0$.

It follows from Lemma 2.7 that the solution $x(t)$ of (1.1) under hypotheses (H1-5) corresponding to arbitrary $d \in M_D$ with arbitrary initial condition $T_r(t_0)x = x_0 \in C^0([-r,0];\Re^n)$ satisfies (A28) and (A29) for all $t \geq t_0$. Therefore system (1.1) is RFC and estimate (3.1) is a direct consequence of (A28) and (3.2) (or (3.7)), as in the previous case.

Proof of (e) $\Rightarrow$ (b):

Let arbitrary $(t_0,x_0) \in \Re^+ \times C^0([-r,0];\Re^n)$ and $d \in \tilde{M}_D$ and consider the solution $x(t)$ of (1.1) with initial condition $T_r(t_0)x = x_0 \in C^0([-r,0];\Re^n)$ corresponding to $d \in \tilde{M}_D$ and defined on $[t_0 - r, +\infty)$. Setting $x(t) := x(t_0 - r)$ for $t \in [t_0 - r - \tau, t_0 - r]$, we may assume that for each time $t \in [t_0, +\infty)$ the unique solution of (1.1) belongs to $C^0([t_0 - r - \tau, t];\Re^n)$. Moreover, we have $\|T_{r+\tau}(t_0)x\|_{r+\tau} = \|T_r(t_0)x\|_r = \|x_0\|_r$. Since (1.1) is Robustly Forward Complete (RFC), by virtue of Lemma 3.5 in [19], there exist functions $\mu \in K^+$, $a \in K_\infty$ and a constant $M \geq 0$, such that for every $(t_0,x_0,d) \in \Re^+ \times C^0([-r,0];\Re^n) \times M_D$, estimate (A2) holds. Without loss of generality we may assume that $\mu \in K^+$ is non-decreasing, so that the following estimate holds:

$$\|T_{r+\tau}(t)x\|_{r+\tau} \leq \mu(t) a(\|x_0\|_r + M), \forall t \geq t_0 \quad (A30)$$

Let $V(t) := V(t,T_{r+\tau}(t)x)$, which is a lower semi-continuous function on $[t_0,+\infty)$. Notice that, by virtue of Lemma 2.4, we obtain:

$$D^+V(t) \leq V^0(t,T_{r+\tau}(t)x;f(t,T_r(0)T_{r+\tau}(t)x,d(t))), \text{ for all } t \geq t_0 \quad (A31)$$

where $D^+V(t) := \limsup_{h \to 0^+} \frac{V(t+h,T_r(t+h)x) - V(t,T_r(t)x)}{h}$. It follows from definition (3.6) that:

$$\text{If } t \geq t_0 + \tau \text{ then } T_{r+\tau}(t)x \in S(t) \quad (A32)$$

Inequality (A31) in conjunction with (A32) and inequality (3.5) gives:

$$D^+V(t) \leq -\gamma(t)\rho(V(t)) + \gamma(t)\mu\left(\int_0^t \gamma(s)ds\right), \text{ for all } t \geq t_0 + \tau \quad (A33)$$

Lemma 2.8 in [22], in conjunction with (A33) and Lemma 5.2 in [17] imply that there exist a function $\sigma(\cdot) \in KL$ and a constant $R > 0$ such that the following inequality is satisfied:

$$V(t) \leq \sigma\left(V(t_0+\tau) + R, \int_{t_0+\tau}^t \gamma(s)ds\right), \forall t \geq t_0 + \tau \quad (A34)$$

It follows from (3.4), (A30) and (A34) that the following estimate holds:

$$V(t) \leq \sigma\left(a_2\left(\beta(t_0+\tau)\mu(t_0+\tau)a(\|x_0\|_r + M)\right) + R, \int_{t_0+\tau}^t \gamma(s)ds\right), \forall t \geq t_0 + \tau \quad (A35)$$

Next, we distinguish the following cases:

1) If (3.4) holds, then (A35) in conjunction with (3.4) and Lemma 3.3 in [19] shows that (1.1) is non-uniformly RGAOS with disturbances $d \in \tilde{M}_D$.

2) If (3.8) holds, then (A35) implies the following estimate:

$$|h(t,x(t))| \leq a_1^{-1}\left(\sigma\left(a_2\left(\beta(t_0+\tau)\mu(t_0+\tau)a(\|x_0\|_r + M)\right) + R, \int_{t_0+\tau}^t \gamma(s)ds\right)\right), \forall t \geq t_0 + \tau$$

and consequently

$$\sup_{\theta \in [-r,0]} |h(t+\theta, x(t+\theta))| \leq a_1^{-1}\left(\sigma\left(a_2\left(\beta(t_0+\tau)\mu(t_0+\tau)a(\|x_0\|_r + M)\right) + R, \int_{t_0+\tau}^{t-r} \gamma(s)ds\right)\right), \forall t \geq t_0 + \tau + r$$

The above estimate, in conjunction with the fact that $H: \Re^+ \times C^0([-r,0]; \Re^n) \to Y$ is equivalent to the finite-dimensional mapping $h$ and Lemma 3.3 in [19] shows that (1.1) is non-uniformly RGAOS with disturbances $d \in \tilde{M}_D$.

The proof is complete. ◁

**Proof of Theorem 3.6:** Implications (a)$\Rightarrow$(b), (d)$\Rightarrow$(c), (c)$\Rightarrow$(e) are obvious. Thus we are left with the proof of implications (b)$\Rightarrow$(d), (c)$\Rightarrow$(a) and (e)$\Rightarrow$(b). The proof of implication of (b)$\Rightarrow$(d) is exactly the same with that of Theorem 3.5 for the special case of the constant function $\beta(t) \equiv 1$. Moreover, the fact that $V$ is $T-periodic$ (or time-independent) if (1.1) is $T-periodic$ (or autonomous) can be shown in the same way as in [22]. The proof of implication (c)$\Rightarrow$(a) is exactly the same with the proof of implication (c)$\Rightarrow$(a) of Theorem 3.5 with the only difference that since $h:[-r,+\infty) \times \Re^n \to \Re^p$ is continuous and $T-periodic$ with $h(t,0) = 0$ for all $t \geq -r$, it follows from Lemma 3.2 in [19] implies that there exist a function $\zeta \in K_\infty$ such that:

$$|h(t-r,x)| \leq \zeta(|x|), \forall (t,x) \in \Re^+ \times \Re^n$$

Finally, the proof of implication (e)$\Rightarrow$(b) follows the same arguments as the proof of implication (e)$\Rightarrow$(b) of Theorem 3.5, with the difference that, by virtue of inequalities (3.12a,b), the function $V(t) := V(t, T_{r+\tau}(t)x)$ satisfies the following differential inequalities:

$$D^+V(t) \leq \beta V(t), \text{ for all } t \geq t_0 \quad (A36)$$

$$D^+V(t) \leq -\rho(V(t)), \text{ for all } t \geq t_0 + \tau \tag{A37}$$

Thus Lemma 4.4 in [28] in conjunction with (A37) shows that there exists $\sigma \in KL$ such that the following inequality is satisfied:

$$V(t) \leq \sigma(V(t_0 + \tau), t - t_0 - \tau), \forall t \geq t_0 + \tau \tag{A38}$$

Moreover, differential inequality (A36) implies $V(t) \leq \exp(\beta(t - t_0))V(t_0)$ for all $t \geq t_0$. Combining the previous estimate with (A38) we obtain:

$$V(t) \leq \omega(V(t_0), t - t_0), \forall t \geq t_0 \tag{A39}$$

where $\omega(s, t) := \max\{\exp(\beta\tau)s, \sigma(s, 0)\}$ for $t \in [0, r)$ and $\omega(s, t) := \max\{\exp(r - t)\exp(\beta\tau)s, \sigma(s, t - r)\}$ for $t \geq r$. From this point the proof can be continued in exactly the same way as in the proof of Theorem 3.5. The proof is complete. ◁

**Proof of Theorem 4.3:**

(a) $\Rightarrow$ (b) : Suppose that there exist functions $\sigma \in KL$, $\beta, \phi \in K^+$, $\rho \in K_\infty$ such that the estimate (4.4) holds for all $(d, u) \in M_D \times M_U$, $(t_0, x_0) \in \Re^+ \times C^0([-r, 0]; \Re^n)$ and $t \geq t_0$. By invoking Lemma 2.3 in [17], there exist functions $p \in K_\infty$ and $\delta \in K^+$ such that $\beta(t)\rho(\phi(t)|u|) \leq p(\delta(t)|u|)$ for all $(t, u) \in \Re^+ \times \Re^m$ and if we set $\gamma(s) := \sigma(p(s), 0) + s$ (that obviously is of class $K_\infty$), the desired (4.1) is a consequence of (4.4) and the previous inequality. Thus statement (b) holds if (1.2) is RFC from the input $u \in M_U$. If the hypothesis that (1.2) is RFC from the input $u \in M_U$ is not included in statement (a) then there exist functions $a \in K_\infty$, $\mu \in K^+$ and a constant $R \geq 0$ such that $a(\mu(t)|x(0)|) \leq \|H(t, x)\|_Y + R$ for all $(t, x) \in \Re^+ \times C^0([-r, 0]; \Re^n)$. It follows from (4.4) and previous definitions that for every $(t_0, x_0) \in \Re^+ \times C^0([-r, 0]; \Re^n)$, $(d, u) \in M_D \times M_U$, the corresponding solution $x(t)$ of (1.2) with $T_r(t_0)x = x_0$ satisfies the following estimate for all $t \geq t_0$:

$$a(\mu(t)|x(t)|) \leq R + \max\left\{\sigma(\beta(t_0)\|x_0\|_r, t - t_0), \sup_{t_0 \leq \tau \leq t} \gamma(\delta(\tau)|u(\tau)|)\right\}$$

The above estimate in conjunction with Definition 2.2 implies that (1.2) is RFC from the input $u \in M_U$.

(b) $\Rightarrow$ (c) : Since (1.2) is RFC from the input $u \in M_U$, by virtue of Lemma 3.5 in [19], there exist functions $q \in K^+$, $a \in K_\infty$ and a constant $R \geq 0$ such that the following estimate holds for all $u \in M_U$ and $(t_0, x_0, d) \in \Re^+ \times C^0([-r, 0]; \Re^n) \times M_D$:

$$\|T_r(t)x\|_r \leq q(t) a\left(R + \|x_0\|_r + \sup_{\tau \in [t_0, t]} |u(\tau)|\right), \forall t \geq t_0 \tag{A40}$$

Using Corollary 10 and Remark 11 in [38], we obtain $\kappa \in K_\infty$ such that $a(rs) \leq \kappa(r)\kappa(s)$ for all $(r, s) \in (\Re^+)^2$. Let $\theta \in K_\infty$ be a locally Lipschitz function that satisfies $\theta(s) \leq \min\{\kappa^{-1}(s); s\}$ for all $s \geq 0$. Moreover, let $\phi(t) := \dfrac{4}{\kappa^{-1}\left(\dfrac{\exp(-t)}{2q(t)}\right)} + \exp(t)q(t)\delta(t)$, $\mu(t) := \dfrac{\exp(-t)}{q(t)}$, where $\delta \in K^+$ is the function involved in (4.1). The previous definitions guarantee that:

$$\text{if } \phi(t)|u| \leq \theta(\|x\|_r) \text{ then } a(4|u|) \leq \frac{1}{2}\mu(t)\|x\|_r \text{ and } \gamma(\delta(t)|u|) \leq \gamma(\mu(t)\|x\|_r) \tag{A41}$$

By virtue of (A40), (4.1) and (A41) it follows that the solution $x(\cdot)$ of (1.2) satisfies the following implication:

$$\phi(\tau)|u(\tau)| \leq \theta(\|T_r(\tau)x\|_r), \text{ a.e. in } [t_0, t] \Rightarrow$$
$$\mu(t)\|T_r(t)x\|_r \leq \exp(-t)a(2R) + \exp(-t)a(4\|x_0\|_r) + \frac{1}{2}\exp(-t)\sup_{t_0 \leq \tau \leq t}\mu(\tau)\|T_r(\tau)x\|_r \quad (A42)$$

$$\phi(\tau)|u(\tau)| \leq \theta(\|T_r(\tau)x\|_r), \text{ a.e. in } [t_0, t] \Rightarrow$$
$$\|H(t, T_r(t)x)\|_Y \leq \sigma(\beta(t_0)\|x_0\|_r, t - t_0) + \sup_{t_0 \leq \tau \leq t}\gamma(\mu(\tau)\|T_r(\tau)x\|_r) \quad (A43)$$

Clearly, (A42) implies:

$$\mu(\tau)\|T_r(\tau)x\|_r \leq \sup_{t_0 \leq s \leq t}\mu(s)\|T_r(s)x\|_r \leq 2a(2R) + 2a(4\|x_0\|_r), \ \forall \tau \in [t_0, t],$$
$$\text{provided that } \phi(\tau)|u(\tau)| \leq \theta(\|T_r(\tau)x\|_r), \text{ a.e. in } [t_0, t] \quad (A44)$$

Notice that every solution $x(\cdot)$ of (4.5) corresponding to some $(d', d) \in \tilde{M}_\Delta$ coincides with the solution of (1.2) corresponding to $u(\cdot) = \dfrac{\theta(\|T_r(\cdot)x\|_r)}{\phi(\cdot)}d'(\cdot)$ initiated from same initial $x_0 \in C^0([-r, 0]; \Re^n)$ at time $t_0 \geq 0$ and corresponding to same $d \in \tilde{M}_D$. Thus, by taking into account (A42), (A43) and (A44), it follows that the solution $x(\cdot)$ of (4.5) satisfies:

$$\mu(t)\|T_r(t)x\|_r \leq \exp(-t)a(2R) + \exp(-t)a(4\|x_0\|_r) + \frac{1}{2}\exp(-t)\sup_{t_0 \leq \tau \leq t}\mu(\tau)\|T_r(\tau)x\|_r,$$
$$\forall t \geq t_0, \ (d', d) \in \tilde{M}_\Delta, \ (t_0, x_0) \in \Re^+ \times C^0([-r, 0]; \Re^n) \quad (A45)$$

$$\|H(t, T_r(t)x)\|_Y \leq \sigma(\beta(t_0)\|x_0\|_r, t - t_0) + \sup_{t_0 \leq \tau \leq t}\gamma(\mu(\tau)\|T_r(\tau)x\|_r),$$
$$\forall t \geq t_0, \ (d', d) \in \tilde{M}_\Delta, \ (t_0, x_0) \in \Re^+ \times C^0([-r, 0]; \Re^n) \quad (A46)$$

$$\mu(t)\|T_r(t)x\|_r \leq 2a(2R) + 2a(4\|x_0\|_r), \ \forall t \geq t_0, \ (d', d) \in \tilde{M}_\Delta, \ (t_0, x_0) \in \Re^+ \times C^0([-r, 0]; \Re^n) \quad (A47)$$

Consider the functions $c(h, T, s) := \sup\{\mu(t_0 + h)\|T_r(t_0 + h)x\|_r \ ; \ (d', d) \in \tilde{M}_\Delta, \|x_0\|_r \leq s, t_0 \in [0, T]\}$ and $b(h, T, s) := \sup\{\|H(t_0 + h, T_r(t_0 + h)x)\|_Y \ ; \ (d', d) \in \tilde{M}_\Delta, \|x_0\|_r \leq s, t_0 \in [0, T]\}$, where $x(\cdot)$ denotes the solution of (4.5) corresponding to some $(d', d) \in \tilde{M}_\Delta$. Next we show that $\lim_{h \to +\infty} b(h, T, s) = \lim_{h \to +\infty} c(h, T, s) = 0$, for all $(T, s) \in (\Re^+)^2$. Clearly, by (A46), (A47) and definitions of $c, b$ we have

$$c(t, T, s) \leq 2a(2R) + 2a(4s) \ ; \ b(t, T, s) \leq \sigma(\max_{0 \leq \tau \leq T}\beta(\tau)s, 0) + \gamma(2a(2R) + 2a(4s)), \ \forall t \geq 0$$

By virtue of the above estimates the mappings $c, b$ are bounded for each fixed $(T, s) \in (\Re^+)^2$ and thus the limits $\limsup_{h \to +\infty} c(h, T, s) = \rho$ and $\limsup_{h \to +\infty} b(h, T, s) = l$ are well defined and finite. We show that $\rho = l = 0$. Indeed, for every $\varepsilon > 0$ there exists $\tau = \tau(\varepsilon, T, s) \geq 0$ such that

$$c(h, T, s) \leq \rho + \varepsilon, \ \forall h \geq \tau \quad (A48)$$

Again recall definitions of $c, b$ above and (A45), (A46), (A47) in conjunction with (A48), which imply

$$b(h,T,s) \leq \sigma\left(\exp(T+\tau)\max_{0\leq\xi\leq T+\tau}\beta(\xi)\max_{0\leq\xi\leq T+\tau}q(\xi)[2a(2R)+2a(4s)], h-\tau\right)+\gamma(\rho+\varepsilon), \text{ for all } h \geq \tau$$

$$c(h,T,s) \leq \exp(-h)a(2R)+\exp(-h)a\left(4\exp(T+\tau)\max_{0\leq\xi\leq T+\tau}q(\xi)[2a(2R)+2a(4s)]\right)+\frac{1}{2}\exp(-h)(\rho+\varepsilon), \text{ for all } h \geq \tau$$

Clearly, the above inequalities imply that $\rho = 0$ as well as $l \leq \zeta(\varepsilon)$ for all $\varepsilon > 0$. Consequently, we must have $\limsup_{h\to+\infty} b(h,T,s) = l = 0$.

The fact that $\lim_{h\to+\infty} b(h,T,s) = \lim_{h\to+\infty} c(h,T,s) = 0$ in conjunction with Lemma 3.3 in [19] shows that (4.5) is non-uniformly in time RGAOS with disturbances $(d',d) \in \tilde{M}_\Delta$.

(c) $\Rightarrow$ (d) : Suppose that (4.5) is non-uniformly in time RGAOS with disturbances $(d',d) \in \tilde{M}_\Delta$. Theorem 3.5 (statement (d)) implies that there exists a continuous mapping $(t,x) \in \Re^+ \times C^0([-r,0];\Re^n) \to V(t,x) \in \Re^+$, which is almost Lipschitz on bounded sets, with the following properties:

-- there exist functions $a_1, a_2 \in K_\infty$, $\beta \in K^+$ such that:

$$a_1(\|H(t,x)\|_Y + \mu(t)\|x\|_r) \leq V(t,x) \leq a_2(\beta(t)\|x\|_r), \quad \forall (t,x) \in \Re^+ \times C^0([-r,0];\Re^n) \tag{A49}$$

-- it holds that:

$$V^0\left(t,x; f\left(t,x,\frac{\theta(\|x\|_r)}{\phi(t)}d',d\right)\right) \leq -V(t,x), \quad \forall (t,x,d',d) \in \Re^+ \times C^0([-r,0];\Re^n) \times \Delta \tag{A50}$$

Notice that inequality (A50) implies the following inequality:

$$V^0(t,x; f(t,x,u,d)) \leq -V(t,x),$$
$$\text{for all } (t,x,u,d) \in \Re^+ \times C^0([-r,0];\Re^n) \times U \times D \text{ with } \phi(t)|u| \leq \theta(\|x\|_r) \tag{A51}$$

Using property (P1) of Definition 2.5 for the continuous mapping $(t,x) \in \Re^+ \times C^0([-r,0];\Re^n) \to V(t,x) \in \Re^+$, we obtain for all $(t,x,u,d) \in \Re^+ \times C^0([-r,0];\Re^n) \times U \times D$:

$$\left|V^0(t,x; f(t,x,u,d)) - V^0(t,x; f(t,x,0,d))\right| \leq M(t+\|x\|_r+1)|f(t,x,u,d) - f(t,x,0,d)|$$

The above inequality in conjunction with (2.6) implies that the following inequality holds for all $(t,x,u,d) \in \Re^+ \times C^0([-r,0];\Re^n) \times U \times D$:

$$\left|V^0(t,x; f(t,x,u,d)) - V^0(t,x; f(t,x,0,d))\right| \leq M(t+\|x\|_r+1)L_U(t,\|x\|_r+|u|)|u| \tag{A52}$$

Define

$$\psi(t,s) := \sup\left\{M(t+\|x\|_r+1)L_U(t,\|x\|_r+|u|)|u| \,;\, \|x\|_r \leq \theta^{-1}(\phi(t)s), |u| \leq s\right\} \tag{A53}$$

Without loss of generality we may assume that the function $\phi \in K^+$ is non-decreasing. Clearly, $\psi: \Re^+ \times \Re^+ \to \Re^+$ is a mapping with $\psi(t,0) = 0$ for all $t \geq 0$, such that (i) for each fixed $t \geq 0$, the mapping $\psi(t,\cdot)$ is non-decreasing;

(ii) for each fixed $s \geq 0$, the mapping $\psi(\cdot, s)$ is non-decreasing and (iii) $\lim_{s \to 0^+} a(t,s) = 0$, for all $t \geq 0$. Hence, by employing Lemma 2.3 in [17], we obtain functions $a_3 \in K_\infty$ and $\delta \in K^+$ such that $\psi(t,s) \leq a_3(\delta(t)s)$.

We next establish inequality (4.7), with $a_3$ as previously, by considering the following two cases:

* $\theta^{-1}(\phi(t)|u|) \leq \|x\|_r$. In this case inequality (4.7) is a direct consequence of (A51).
* $\theta^{-1}(\phi(t)|u|) \geq \|x\|_r$. In this case, by virtue of inequalities (A51), (A52), definition (A53) and definition of $a_3$, we have: $V^0(t,x; f(t,x,u,d)) \leq V^0(t,x; f(t,x,0,d)) + \psi(t,|u|) \leq -V(t,x) + a_3(\delta(t)|u|)$, which implies (4.7).

(d) $\Rightarrow$ (e) : Notice that (4.7) implies (4.9) with $a(s) := 2a_3(s)$ and $\rho(s) := \frac{1}{2}s$.

(e) $\Rightarrow$ (a) : Theorem 4.6 implies that system (1.2) is RFC from the input $u \in M_U$ and that (4.20) holds. Next, we distinguish the following cases:

1) If (4.8) holds, then (4.4) is a direct consequence of (4.20) and (4.8).
2) If (4.10) holds, then (4.20) implies the following estimate:

$$|h(t,x(t))| \leq \max\left\{ a_1^{-1}\left(\sigma\left(a_2(\beta(t_0)\|x_0\|_r), t-t_0\right)\right), \sup_{t_0 \leq s \leq t} a_1^{-1}\left(\sigma\left(\zeta(\delta(s)|u(s)|), t-s\right)\right) \right\}, \forall t \geq t_0$$

Since $h: [-r, +\infty) \times \Re^n \to \Re^p$ is continuous with $h(t,0) = 0$ for all $t \geq -r$, it follows from Lemma 3.2 in [19] implies that there exist functions $p \in K_\infty$ and $\phi \in K^+$ such that:

$$|h(t-r, x)| \leq p(\phi(t)|x|), \quad \forall (t,x) \in \Re^+ \times \Re^n$$

Combining the two previous inequalities we obtain:

$$\sup_{\theta \in [-r,0]} |h(t+\theta, x(t+\theta))| \leq \max\left\{ \omega\left(q(t_0)\|x_0\|_r, t-t_0\right), \sup_{t_0 \leq s \leq t} \omega\left(q(s)\zeta(\delta(s)|u(s)|), t-s\right) \right\}, \forall t \geq t_0$$

where $q(t) := 1 + \beta(t) + \max_{t \leq \tau \leq t+r} \phi(\tau)$ $\omega(s,t) := \max\left\{ p(s), a_1^{-1}(\sigma(s+a_2(s), 0)) \right\}$ for $t \in [0, r)$ and $\omega(s,t) := \max\left\{ \exp(r-t)p(s), a_1^{-1}(\sigma(s+a_2(s), t-r)) \right\}$ for $t \geq r$. The above estimate, in conjunction with the fact that $H: \Re^+ \times C^0([-r,0]; \Re^n) \to Y$ is equivalent to the finite-dimensional mapping $h$ shows that (1.2) satisfies inequality (4.4).

The proof is complete. ◁

**Proof of Theorem 4.4:** The proof of implications (a) $\Rightarrow$ (b), (d) $\Rightarrow$ (e) and (e) $\Rightarrow$ (a) follow the same methodology as in the proof of Theorem 4.3. Particularly, in the proof of implication (e) $\Rightarrow$ (a), we use in addition the fact that since $h: [-r, +\infty) \times \Re^n \to \Re^p$ is continuous and $T-periodic$ with $h(t,0) = 0$ for all $t \geq -r$, it follows from Lemma 3.2 in [19] implies that there exist a function $p \in K_\infty$ such that:

$$|h(t-r, x)| \leq p(|x|), \quad \forall (t,x) \in \Re^+ \times \Re^n$$

The proof of implication (c) $\Rightarrow$ (d) differs from the corresponding proof in Theorem 4.3 in the definition of $\psi$. Specifically, we first notice that the fact that $V$ is $T-periodic$, implies that $V^0(t,x; f(t,x,u,d))$ is $T-periodic$. Using property (P1) of Definition 2.5 for the continuous mapping $(t,x) \in \Re^+ \times C^0([-r,0]; \Re^n) \to V(t,x) \in \Re^+$, we obtain for all $(t,x,u,d) \in \Re^+ \times C^0([-r,0]; \Re^n) \times U \times D$:

$$\left|V^0(t,x;f(t,x,u,d))-V^0(t,x;f(t,x,0,d))\right| \leq M\left(T+\|x\|_r+1\right)\left|f(t,x,u,d)-f(t,x,0,d)\right|$$

The above inequality in conjunction with (2.6) and the fact that $f$ is $T-periodic$ implies that the following inequality holds for all $(t,x,u,d) \in \Re^+ \times C^0([-r,0];\Re^n) \times U \times D$:

$$\left|V^0(t,x;f(t,x,u,d))-V^0(t,x;f(t,x,0,d))\right| \leq M\left(T+\|x\|_r+1\right)L_U(T,\|x\|_r+|u|)|u|$$

We next define:

$$\psi(s) := \sup\left\{M\left(T+\|x\|_r+1\right)L_U(T,\|x\|_r+|u|)|u|\,;\,\|H(t,x)\|_Y \leq \theta^{-1}(s)\,,|u| \leq s\right\}$$

Notice that hypothesis (S8) implies that $\psi(s) \leq a_3(s)$ for all $s \geq 0$, where $a_3(s) := \tilde{M}\left(T+R+a(\theta^{-1}(s))\right)L_U(T,R+a(\theta^{-1}(s))+s)s$ and $\tilde{M}(s)$ is a continuous positive function which satisfies $\tilde{M}(s) \geq M(s)$ for all $s \geq 0$. From this point the proof of implication (c) $\Rightarrow$ (d) is exactly the same as in Theorem 4.3 (i.e., by distinguishing the cases $\theta^{-1}(|u|) \leq \|H(t,x)\|_Y$ and $\theta^{-1}(|u|) \geq \|H(t,x)\|_Y$).

Finally, we continue with the proof of implication (b) $\Rightarrow$ (c). Without loss of generality we may assume that $\gamma \in K_\infty$. Let $\theta \in K_\infty$ be a locally Lipschitz function that satisfies $\theta(s) \leq \gamma^{-1}\left(\frac{1}{2}s\right)$ for all $s \geq 0$. By virtue of (4.1) and hypothesis (S8) it follows that the solution $x(\cdot)$ of (1.2) satisfies the following implications:

$$|u(\tau)| \leq \theta\left(\|H(\tau,T_r(\tau)x)\|_Y\right), \text{ a.e. in } [t_0,t] \Rightarrow$$
$$\|H(t,T_r(t)x)\|_Y \leq \max\left\{\sigma(\|x_0\|_r,t-t_0);\frac{1}{2}\sup_{t_0 \leq \tau \leq t}\|H(\tau,T_r(\tau)x)\|_Y\right\} \tag{A54}$$

$$|u(\tau)| \leq \theta\left(\|H(\tau,T_r(\tau)x)\|_Y\right), \text{ a.e. in } [t_0,t] \Rightarrow$$
$$\|T_r(t)x\|_r \leq R+a\left(\sigma(\|x_0\|_r,t-t_0)+\frac{1}{2}\sup_{t_0 \leq \tau \leq t}\|H(\tau,T_r(\tau)x)\|_Y\right) \tag{A55}$$

Proceeding in exactly the same way as in the proof of Theorem 4.3, it can be that for all $(d',d) \in M_\Delta$, $(t_0,x_0) \in \Re^+ \times C^0([-r,0];\Re^n)$, the corresponding solution $x(\cdot)$ of (4.12) satisfies the following estimates for all $t \geq t_0$:

$$\|H(t,T_r(t)x)\|_Y \leq \max\left\{\sigma(\|x_0\|_r,t-t_0);\frac{1}{2}\sup_{t_0 \leq \tau \leq t}\|H(\tau,T_r(\tau)x)\|_Y\right\} \tag{A56}$$

$$\|T_r(t)x\|_r \leq R+a\left(\sigma(\|x_0\|_r,t-t_0)+\frac{1}{2}\sup_{t_0 \leq \tau \leq t}\|H(\tau,T_r(\tau)x)\|_Y\right) \tag{A57}$$

$$\|H(t,T_r(t)x)\|_Y \leq \sigma(\|x_0\|_r,0) \text{ and } \|T_r(t)x\|_r \leq R+a(2\sigma(\|x_0\|_r,0)) \tag{A58}$$

Consequently, (A58) and Definition 2.1 imply that system (4.12) is RFC. Moreover, using estimates (A56) and (A58), it follows that for all $t \geq t_0$ we have:

$$\|H(t,T_r(t)x)\|_Y \leq \inf_{t_0 \leq \xi \leq t}\max\left\{\sigma(R+a(2\sigma(\|x_0\|_r,0)),t-\xi);\frac{1}{2}\sup_{\xi \leq \tau \leq t}\|H(\tau,T_r(\tau)x)\|_Y\right\} \tag{A59}$$

Lemma 2.10 in conjunction with inequality (A59) guarantees the existence of $\omega \in KL$ such that for all $t \geq t_0$ we have:

$$\|H(t,T_r(t)x)\|_Y \le \omega\!\left(R + a\!\left(2\sigma(\|x_0\|_r,0)\right), t - t_0\right) \tag{A60}$$

Combining estimate (A58) with (A60) we obtain for all $t \ge t_0$:

$$\|H(t,T_r(t)x)\|_Y \le \kappa(\|x_0\|_r, t - t_0) \tag{A61}$$

where $\kappa(s,t) := (\sigma(s,0))^{\frac{1}{2}} (\omega(R + a(2\sigma(s,0)), t))^{\frac{1}{2}}$ (notice that $\kappa$ is of class $KL$). Estimate (A61) in conjunction with Definition 3.3 and the fact that system (4.12) is RFC shows that (4.12) is URGAOS with disturbances $(d',d) \in \tilde{M}_\Delta$. The proof is complete. ◁

**Proof of Theorem 4.6:** Consider a solution of (1.2) under hypotheses (S1-7) corresponding to arbitrary $(u,d) \in M_U \times M_D$ with initial condition $T_r(t_0)x = x_0 \in C^1([-r,0];\mathfrak{R}^n)$. By virtue of Lemma 2.6, for every $T \in (t_0, t_{\max})$, the mapping $[t_0,T] \ni t \to V(t,T_r(t)x)$ is absolutely continuous. It follows from (4.19) and Lemma 2.4 that there exists a set $I \subset [t_0,T]$ of zero Lebesgue measure such that the following implication holds for all $t \in [t_0,T] \setminus I$:

$$V(t,T_r(t)x) \ge \zeta(\delta(t)|u(t)|) \Rightarrow \frac{d}{dt}(V(t,T_r(t)x)) \le -\rho(V(t,T_r(t)x))$$

Lemma 4.7 implies the existence of a continuous function $\sigma$ of class $KL$, with $\sigma(s,0) = s$ for all $s \ge 0$ such that:

$$V(t,T_r(t)x) \le \max\!\left\{ \sigma(V(t_0,T_r(t_0)x), t - t_0), \sup_{t_0 \le s \le t} \sigma\!\left(\zeta(\delta(s)|u(s)|), t - s\right) \right\}, \;\forall t \in [t_0,T] \tag{A62}$$

Notice that for the case that (1.2) is RFC from the input $u \in M_U$ then $t_{\max} = +\infty$. For the case that there exist functions $p \in K_\infty$, $\mu \in K^+$ and a constant $R \ge 0$ such that $p(\mu(t)|x(0)|) \le V(t,x) + R$ for all $(t,x) \in \mathfrak{R}^+ \times C^0([-r,0];\mathfrak{R}^n)$, combining the previous inequality and (A62) we obtain for every $T \in (t_0, t_{\max})$:

$$p(\mu(t)|x(t)|) \le R + \max\!\left\{ \sigma(V(t_0,T_r(t_0)x), t - t_0), \sup_{t_0 \le s \le t} \sigma\!\left(\zeta(\delta(s)|u(s)|), t - s\right) \right\}, \;\forall t \in [t_0,T] \tag{A63}$$

It follows from estimate (A63) that $t_{\max} = +\infty$. Immediate consequence of Lemma 2.7 is that estimate (A62) holds for all $(t_0,x_0,d,u) \in \mathfrak{R}^+ \times C^0([-r,0];\mathfrak{R}^n) \times M_D \times M_U$ and $t \ge t_0$. Moreover, if there exist functions $p \in K_\infty$, $\mu \in K^+$ and a constant $R \ge 0$ such that $p(\mu(t)|x(0)|) \le V(t,x) + R$ for all $(t,x) \in \mathfrak{R}^+ \times C^0([-r,0];\mathfrak{R}^n)$, then Lemma 2.7 implies that estimate (A63) also holds for all $(t_0,x_0,d,u) \in \mathfrak{R}^+ \times C^0([-r,0];\mathfrak{R}^n) \times M_D \times M_U$ and $t \ge t_0$. In this case the fact that system (1.2) is RFC from the input $u \in M_U$ is an immediate consequence of (A63) and Definition 2.2.

Notice that (4.20) is an immediate consequence of (A62) and (4.18). Finally, (i) and (ii) are immediate consequences of (4.20). The proof is complete. ◁

**Proof of Lemma 4.7:** Notice that by virtue of Lemma 4.4 in [28], for each positive definite continuous function $\rho: \mathfrak{R}^+ \to \mathfrak{R}^+$ there exists a continuous function $\sigma$ of class $KL$, with $\sigma(s,0) = s$ for all $s \ge 0$ with the following property: if $y:[t_0,t_1] \to \mathfrak{R}^+$ is an absolutely continuous function and $I \subset [t_0,t_1]$ a set of Lebesgue measure zero such that $\dot{y}(t)$ is defined on $[t_0,t_1] \setminus I$ and such that the following differential inequality holds for all $t \in [t_0,t_1] \setminus I$:

$$\dot{y}(t) \le -\rho(y(t)) \tag{A64}$$

then the following estimate holds for all $t \in [t_0, t_1]$:

$$y(t) \leq \sigma\left(y(t_0), t - t_0\right) \tag{A65}$$

Actually, the statement of Lemma 4.4 in [28] does not guarantee that $\sigma$ is continuous or that $\sigma(s,0) = s$ for all $s \geq 0$, but a close look at the proof of Lemma 4.4 in [28] shows that this is the case when $\rho : \Re^+ \to \Re^+$ is a positive definite continuous function. Moreover, notice that we may continuously extend $\sigma$ by defining $\sigma(s,t) := s \exp(-t)$ for $t < 0$.

Clearly, (4.22) holds for $t = t_0$ and $\sigma$ the function involved in (A65). We next show that (4.22) holds for arbitrary $t \in (t_0, t_1]$.

Let arbitrary $t \in (t_0, t_1]$ and define the functions $\tilde{u}(\tau) = \begin{cases} u(\tau) \text{ for } \tau \in [t_0, t] \\ 0 \text{ if otherwise} \end{cases}$, $\bar{u}(\tau) := \limsup_{\xi \to \tau} \tilde{u}(\xi)$. Notice that $\bar{u}$ is upper semi-continuous on $\tau \in [t_0, t]$ and consequently the function $p(\tau) := y(\tau) - \bar{u}(\tau)$ is lower semi-continuous on $[t_0, t]$. Next define the set:

$$A := \{\tau \in [t_0, t] : y(\tau) \leq \bar{u}(\tau)\} \tag{A66}$$

We distinguish the following cases:

1) $A = \emptyset$. In this case we have $y(\tau) > \bar{u}(\tau)$ for all $\tau \in [t_0, t]$. Since $\bar{u}(\tau) \geq u(\tau)$ for all $\tau \in [t_0, t]$, the previous inequality in conjunction with (4.21) implies that $\dot{y}(\tau) \leq -\rho(y(\tau))$ for all $\tau \in [t_0, t] \setminus N$. Thus in this case Lemma 4.4 in [28] guarantees that estimate (A65) holds.

2) $A \neq \emptyset$ and $\xi := \sup A < t$. In this case there exists a sequence $\tau_i \leq \xi$ with $\tau_i \to \xi$ and $y(\tau_i) - \bar{u}(\tau_i) \leq 0$. Since the function $p(t) = y(t) - \bar{u}(t)$ is lower semi-continuous, we obtain $p(\xi) = \liminf_{\tau \to \xi} p(\tau) \leq 0$ and consequently $y(\xi) \leq \bar{u}(\xi)$. Moreover, notice that by virtue of definition (A66), (4.21) and since $\bar{u}(\tau) \geq u(\tau)$ for all $\tau \in [t_0, t]$, the differential inequality $\dot{y}(\tau) \leq -\rho(y(\tau))$ holds for all $\tau \in (\xi, t] \setminus N$. Consequently, Lemma 4.4 in [28] implies $y(t) \leq \sigma(y(\tau), t - \tau)$ for all $\tau \in (\xi, t]$. By virtue of continuity of $\sigma$ and $y$ we get

$$y(t) \leq \sigma(y(\xi), t - \xi)$$

which combined $y(\xi) \leq \bar{u}(\xi)$ directly implies

$$y(t) \leq \sigma(\bar{u}(\xi), t - \xi) \leq \sup_{t_0 \leq s \leq t} \sigma(\bar{u}(s), t - s) \tag{A67}$$

3) $A \neq \emptyset$ and $\xi := \sup A = t$. In this case there exists a sequence $\tau_i \leq t$ with $\tau_i \to t$ and $y(\tau_i) - \bar{u}(\tau_i) \leq 0$. Since the function $p(t) = y(t) - \bar{u}(t)$ is lower semi-continuous, we obtain $p(t) = \liminf_{\tau \to t} p(\tau) \leq 0$ and consequently $y(t) \leq \bar{u}(t)$. Moreover, since $\sigma(s,0) = s$ for all $s \geq 0$, it holds that:

$$y(t) \leq u(t) = \sigma(\bar{u}(t), 0) \leq \sup_{t_0 \leq s \leq t} \sigma(\bar{u}(s), t - s)$$

Combining all the above cases, we may conclude that

$$y(t) \leq \max\left\{\sigma(y(t_0), t - t_0), \sup_{t_0 \leq s \leq t} \sigma(\bar{u}(s), t - s)\right\} \tag{A68}$$

Let $M := \sup_{t_0 \leq s \leq t} u(s)$. For each $\varepsilon > 0$, there exists $\delta > 0$ such that $\sigma(s, \tau - \delta) - \sigma(s, \tau) < \varepsilon$ for all $(s, \tau) \in [0, M] \times [0, t]$. Notice that since $\bar{u}(\tau) := \limsup_{\xi \to \tau} \tilde{u}(\xi)$ and $\tilde{u}(\tau) = \begin{cases} u(\tau) \text{ for } \tau \in [t_0, t] \\ 0 \text{ if otherwise} \end{cases}$, it follows that $\bar{u}(s) \leq \sup\{u(r) : \max(s - \delta, t_0) \leq r \leq \min(s + \delta, t)\}$, for all $s \in [t_0, t]$. The previous inequalities imply:

$$\sigma(\bar{u}(s), t - s) \leq \sup\{\sigma(u(r), t - s) : \max(s - \delta, t_0) \leq r \leq \min(s + \delta, t)\}$$
$$\leq \sup\{\sigma(u(r), t - r - \delta) : \max(s - \delta, t_0) \leq r \leq \min(s + \delta, t)\}$$
$$\leq \sup\{\sigma(u(r), t - r) : \max(s - \delta, t_0) \leq r \leq \min(s + \delta, t)\} + \varepsilon$$
$$\leq \sup_{t_0 \leq r \leq t} \sigma(u(r), t - r) + \varepsilon$$

The above inequality in conjunction with (A68) imply that for each $\varepsilon > 0$, it holds that:

$$y(t) \leq \max\left\{\sigma(y(t_0), t - t_0), \sup_{t_0 \leq s \leq t} \sigma(u(s), t - s)\right\} + \varepsilon$$

Since $\varepsilon > 0$ is arbitrary, we conclude that the above estimate directly implies (4.22). The proof is complete. ◁